\documentclass[headings=standardclasses]{scrartcl}

\ifdefined\directlua
\usepackage{fontspec}
\else
\usepackage[T1,T2A]{fontenc}
\fi
\usepackage[russian,UKenglish]{babel}
\usepackage{amssymb,amsmath,amsopn,amsthm,amsfonts,bm}
\usepackage{stmaryrd}
\usepackage[final, pdftitle={}, pdfauthor={Daniel Seibel},
pdfkeywords={}, colorlinks, bookmarksnumbered, bookmarksopen,
bookmarksopenlevel=1, pdfstartview=FitH, linkcolor=black,
citecolor=black, urlcolor=black, filecolor=black ]{hyperref}
\usepackage{csquotes}
\usepackage[backend=bibtex8,style=numeric-comp,sorting=none,
isbn=false,url=false,doi=false,eprint=false]{biblatex}
\setkomafont{labelinglabel}{\bfseries}

\usepackage{graphicx}
\usepackage{subcaption}
\usepackage{color}
\graphicspath{{figures/}}
\DeclareGraphicsExtensions{.pdf, .jpg, .png, .ps, .eps}
\usepackage{booktabs,colortbl}
\usepackage[font={small,it}]{caption}
\usepackage{rotating}
\usepackage{tabularx}

\theoremstyle{plain}

\theoremstyle{definition}
\newtheorem{definition}{Definition}

\newtheorem{remark}{Remark}

\DeclareMathOperator{\supp}{supp}

\DeclareMathOperator{\vspan}{span}

\DeclareMathOperator{\sgn}{sgn}
\DeclareMathOperator{\arcsinh}{arcsinh}

\newcommand{\norm}[1]{\left\lVert#1\right\rVert}
\newcommand{\seminorm}[1]{\left\lvert#1\right\rvert}
\newcommand{\dotprod}[1]{\left\langle#1\right\rangle}

\begin{document}

\title{Almost complete analytical integration in Galerkin BEM}
\author{Daniel Seibel
  \thanks{\emph{Faculty of Mathematics and Computer Science, Campus E1.1,
      Saarland University, 66123 Saarbr{\"u}cken, Germany},
    E-Mail:\ \href{mailto:seibel@num.uni-sb.de}{\texttt{seibel@num.uni-sb.de}}
  }\footnotemark[1]
}

\date{\today}

\maketitle

\begin{abstract}  
  In this work, semi-analytical formulae for the numerical evaluation of surface
  integrals occurring in Galerkin boundary element methods (BEM) in 3D are
  derived. The integrals appear as the entries of BEM matrices and are formed
  over pairs of surface triangles. Since the integrands become singular if the
  triangles have non-empty intersection, the transformation presented
  in~\cite{SauterSchwab2011} is used to remove the singularities. It is shown
  that the resulting integrals admit analytical formulae if the triangles are
  identical or share a common edge. Moreover, the four-dimensional integrals are
  reduced to one- or two-dimensional integrals for triangle pairs with common
  vertices or disjoint triangles respectively. The efficiency and accuracy of
  the formulae is demonstrated in numerical experiments.
\end{abstract}

\section{Introduction}\label{sec:intro}
Whenever unbounded domains appear in the modelling of physical problems,
boundary element methods (BEM) present a particularly effective tool for the 
numerical simulation. Instead of discretising the underlying boundary
value problem directly, BEM operate on the corresponding integral equations
posed on the boundary. Hence, shape functions are defined on the surface and not
in the volume, which results in less degrees of freedoms overall. BEM are
therefore applied in numerous fields of science, ranging from computational
elasticity to electromagnetic and acoustic scattering.

However, the price one pays is the occurrence of dense matrices that are
expensive to calculate. A Galerkin approximation of the Laplace equation with
trial and test functions $\varphi$ and $\psi$ requires the calculation of
integrals over surface triangles $\sigma$ and $\tau$,
\[
  I = \int\limits_{\tau} \int\limits_{\sigma}
  \frac{1}{4\pi \seminorm{\bm{x}-\bm{y}}}
  \, \varphi(\bm{y}) \,dS(\bm{y})\, \psi(\bm{x}) \,dS(\bm{x}),
\]
which do not vanish. In addition, the kernel function is singular
at $\bm{x}=\bm{y}$, so standard quadrature rules for the approximation of $I$
may perform poorly. Whereas the problem of fully populated matrices can be
solved by hierarchical low-rank approximation~\cite{Bebendorf2008}, several
algorithms for the efficient calculation of $I$ are available in the literature.
The general strategy is to use coordinate transformations that render the integral
suitable for numerical integration. Approaches based on polar or Duffy
coordinates~\cite{Duffy1982,MousaviSukumar2010} yield regular integrals for a
wide selection of kernel functions~\cite{SchwabWendland1992}, which can be
approximated by product quadrature rules~\cite{Sauter1998}. However, since the
integrals are essentially four-dimensional, numerical quadrature is
expensive. In order to reduce the computational effort, analytical integration
can be carried out for specific kernels to obtain lower dimensional
integrals~\cite{RjasanowSteinbach2007,Taylor2003}.
If the integral is only defined in weak sense as a finite-part or Hadamard
integral, then integration by parts is often the appropriate
solution~\cite{DautrayLions1990}.

The main contribution of this article is the derivation of analytical formulae
for the complete integration of $I$ for singular cases. Our approach is based on
the regularisation method by S. Erichsen and S.\ A. Sauter~\cite{Sauter1998}
which removes the singularities by applying a variant of the Duffy
transformation. We show that the resulting representation admits closed formulae
of $I$ for identical triangles as well as triangles with a common edge and
reduce it to a one-dimensional integral for triangles with a common vertex.  

\section{Preliminaries}\label{sec:prelim}
We consider the numerical solution of the Laplace problem
\begin{equation}\label{eq:pde}
  \begin{aligned}
    -\Delta u &= 0 && \textup{in } \Omega,\\
    u &= g && \textup{on } \Gamma = \partial \Omega,
  \end{aligned}
\end{equation}
in a domain $\Omega$ with bounded Lipschitz boundary $\Gamma$. If $\Omega$ is
unbounded, we assume the radiation condition
\[
  \seminorm{u(\bm{x})} \in \mathcal{O}\left(\seminorm{\bm{x}}^{-1}\right)\quad
  \textup{for } \seminorm{\bm{x}} \to \infty.
\]
The representation formula expresses the solution $u$ in terms of its boundary
values only,
\[
  u(\bm{x}) = \int_{\Gamma} u^\ast(\bm{x},\bm{y})\,
  \partial_{\bm{n}} u(\bm{y})\, dS(\bm{y})
  - \int_{\Gamma} \partial_{\bm{n}(\bm{y})} u^\ast(\bm{x},\bm{y})\, g(\bm{y})
  \, dS(\bm{y}), \quad \bm{x} \in \Omega.
\]
Here, $\bm{n}$ is the unit normal to $\Gamma$ pointing outwards $\Omega$ and
$u^\ast$ is the fundamental solution of the Laplace operator,
\[
  u^\ast (\bm{x},\bm{y}) = \frac{1}{4\pi \seminorm{\bm{y}-\bm{x}}}.
\]
The boundary value problem is hence reduced to the problem of finding the unknown
Neumann trace $t=\partial_{\bm{n}} u$. To this end, we take the traces in the
representation formula and insert the Dirichlet condition to obtain the boundary
integral equation
\begin{equation}\label{eq:bie}
  \mathcal{V} t = \left(\frac12 \mathcal{I} + \mathcal{K}\right) g\quad
  \textup{on } \Gamma,
\end{equation}
where the layer potentials are defined by
\[
  \begin{aligned}   
    \left(\mathcal{V} w\right)(\bm{x}) &=
    \int_\Gamma u^\ast(\bm{x},\bm{y})\, w(\bm{y}) \, dS(\bm{y}),&
    \left(\mathcal{K} w\right)(\bm{x}) &=
    \int_\Gamma \partial_{\bm{n}(\bm{y})} u^\ast(\bm{x},\bm{y})\,
    w(\bm{y}) \, dS(\bm{y}).
  \end{aligned}
\]
Neumann or mixed boundary conditions can be treated similarly. We refer
to~\cite{McLean2000} for more details.

For the numerical solution of~\eqref{eq:bie} with BEM, we discretise the
boundary with finite elements.

\begin{definition}[Mesh]
  A mesh $\left(\Gamma_h,\mathcal{T}_h\right)$ (or simply $\Gamma_h$) is a
  finite collection of non-empty and open elements $\tau \subset \Gamma_h$ which
  satisfies: 
  \begin{enumerate}
  \item $\mathcal{T}_h = {\left\{ \tau_n \right\}}_{n=1}^N$ is a triangulation of
    $\Gamma_h$, i.e. 
    \[
      \Gamma_h = \bigcup_{n=1}^N \bar{\tau}_n.
    \]
    The intersection $\bar{\tau}_n \cap \bar{\tau}_m$ of two distinct elements
    is either empty or consists of a common vertex or edge.
  \item Each $\tau$ in $\mathcal{T}_h$ is a flat triangle with vertices
    $\bm{p}_1,\bm{p}_2$ and $\bm{p}_3$. The reference mapping
    \[
      \chi_\tau: \pi \to \tau, \quad
      \chi_\tau(x_1,x_2)=\bm{p}_1 +x_1 (\bm{p}_2-\bm{p}_1) + x_2
      (\bm{p}_3-\bm{p}_2),
    \]
    parametrises $\tau$ by the reference triangle 
    \[
      \pi = \left\{ (x_1,x_2)
        \mid 0<x_1<1, 0 < x_2 < x_1
      \right\} \subset \mathbb{R}^2.
    \]    
    We denote by
    \[
      \bm{J}_\tau=\left(\bm{p}_2-\bm{p}_1\mid
        \bm{p}_3-\bm{p}_2\right) \in \mathbb{R}^{3 \times 2},\quad
      g_\tau=\sqrt{\det\left(\bm{J}_\tau^\top \bm{J}_\tau\right)}
    \]
    the Jacobian and the Gram determinant of $\chi_\tau$ respectively and
    assume $g_\tau \neq 0$.
  \end{enumerate}
\end{definition}

\begin{remark}
  Certainly, the particular choice of the reference element $\pi$ is not
  important for our approach. The reason why we use a non-standard $\pi$
  nonetheless lies in the fact that the presentation in Section~\ref{sec:reg}
  becomes simpler and it is in accordance with the literature referenced there.
\end{remark}

On the triangular mesh, we define piece-wise constant and piece-wise linear
ansatz functions. 

\begin{definition}[Boundary element spaces]
  For $p=0,1$, we denote by
  \[
    S_h^0(\pi) = \{1 \},\quad
    S_h^1(\pi) = \{1-x_1,x_1 -x_2,x_2\}
  \]
  the set of reference functions and by
  \[
    S_h^p(\tau) = \vspan \left\{ \varphi \circ \chi_\tau^{-1} :
      \varphi \in S_h^p(\pi)\right\}
  \]
  the local boundary element space on $\tau$. We define the global space by
  gluing the local spaces together, i.e.
  \[
    S_h^p = \left\{\varphi:\Gamma_h\to\mathbb{R} : \varphi_{|\tau}\in
      S_h^p(\tau)\ \forall \tau \in \mathcal{T}_h\right\}.
  \]
  For $p=1$, we moreover require that the functions $\varphi$ are continuous.
\end{definition}

We choose the Lagrangian basis
\[
  \varphi_n^0 (\bm{x}) = \left\{
    \begin{aligned}
      &1 &&\textup{if } \bm{x} \in \tau_n,\\
      &0 &&\textup{else}
    \end{aligned}
    \right\},\quad
  \varphi_m^1(\bm{x}_j) = \left\{
    \begin{aligned}
      &1 &&\textup{if } j = m,\\
      &0 &&\textup{else}
    \end{aligned}
    \right\},
\]
where ${\{\bm{x}_j\}}_{j=1}^M$ denotes the set of vertices in $\Gamma_h$.
Then, the ansatz
\[
  \begin{aligned}
    t_h &= \sum_{n=1}^{N} \bm{t}_n\, \varphi_n^{0}
    \in S_h^{0},& \bm{t} \in\mathbb{R}^{N},&&
    g_h &= \sum_{m=1}^{M} \bm{g}_m\, \varphi_m^1
    \in S_h^1,& \bm{g} \in\mathbb{R}^{M},
  \end{aligned}
\]
for the approximate boundary data leads to the Galerkin approximation
of~\eqref{eq:bie}: Find $\bm{t} \in \mathbb{R}^{N}$ such that
\begin{equation}\label{eq:galapx}
  \bm{V} \bm{t} = \left(\frac 12 \bm{M} + \bm{K} \right) \bm{g},
\end{equation}
where the matrices $\bm{M} \in \mathbb{R}^{N \times M},
\bm{V} \in \mathbb{R}^{N \times N}$ and $\bm{K} \in \mathbb{R}^{N \times M}$ are
given by
\[
  M[n,m] = \dotprod{\varphi_n^{0},\varphi_m^{1}},\quad
  V[n,i] = \dotprod{\varphi_n^{0},\mathcal{V} \varphi_i^{0}},\quad
  K[n,m] = \dotprod{\varphi_n^{0},\mathcal{K} \varphi_m^{1}},
\]
with $i,n=1,\ldots,N$ and $m=1,\ldots,M$. The brackets symbolise the usual
$L_2$-inner product
\[
  \dotprod{u,v} = \int\limits_{\Gamma_h} u(\bm{x})\,v(\bm{x})\,dS(\bm{x}).
\]
The boundary integral equation is now reduced to a system of linear equations,
which can be solved efficiently with direct or iterative methods.

\section{Integral Regularisation}\label{sec:reg}
The entries of $\bm{V}$ an $\bm{K}$ are of the form
\begin{equation}\label{eq:genentry}
  \int\limits_{\Gamma_h\times\Gamma_h} k(\bm{x},\bm{y})\, \varphi(\bm{y}) \,dS(\bm{y})
  \psi(\bm{x}) \,dS(\bm{x})
  = \sum_{\sigma,\tau \in \mathcal{T}_h}
  \int\limits_\sigma \int\limits_\tau k(\bm{x},\bm{y})\, \varphi(\bm{y}) \,dS(\bm{y})
  \psi(\bm{x}) \,dS(\bm{x}),
\end{equation}
where $k = u^\ast, \partial_{\bm{n}(\bm{y})} u^\ast$ is the kernel function and
$\varphi, \psi$ are trial and test functions respectively. Let $I$ be one of the
summands for the non-trivial case $\tau \subset \supp \varphi$ and $\sigma
\subset \supp \psi$. We transform back to the reference element $\pi$,
\[
  \begin{aligned}
    I &= \int\limits_{\sigma} \int\limits_{\tau} k(\bm{x},\bm{y})\, \varphi(\bm{y})
    \,dS(\bm{y}) \psi(\bm{x}) \,dS(\bm{x}) \\
    &= \int\limits_{\pi \times \pi}
    g_\sigma\,g_\tau\,
    k(\chi_\sigma(\bm{x}),\chi_\tau(\bm{y}))\, 
    \varphi(\chi_\sigma(\bm{x}))\, \psi(\chi_\tau(\bm{y}))\,
    d(\bm{y})\,d(\bm{x}),
  \end{aligned}
\]
and abbreviate the integrand by $q$. Since the kernel function
$k(\bm{x},\bm{y})$ is singular at $\bm{x}=\bm{y}$, the integral needs to be
regularised. We distinguish between four different cases: the intersection
$\bar{\sigma} \cap \bar{\tau}$ may consist either of 
\begin{enumerate}
\item the whole element,
\item exactly one edge,
\item exactly one point,
\item be empty.
\end{enumerate}
In the following, we summarise the regularisation introduced
in~\cite{Sauter1998}. For the most part, we adhere to the version
of~\cite[Chapter 5]{SauterSchwab2011}. 

\subsection{Identical elements}\label{subsec:reg_id}
For identical elements $\sigma = \tau$, we substitute
\[
  \bm{z} = \bm{x} - \bm{y}, \quad
  \bm{Z} = (\bm{x} + \bm{y})/2
\]
such that
\[
  I =  \int\limits_{\Pi}
  \int\limits_{\pi_{\bm{z}}}
  q (\bm{Z}+\bm{z}/2,\bm{Z}-\bm{z}/2)
  \, d\bm{Z}\, d\bm{z}
\]
with
\[
  \Pi = \left\{z=\bm{x}-\bm{y} : \bm{x},\bm{y}\in
    \pi \right\},\quad
  \pi_{\bm{z}}=(\pi - \bm{z}/2) \cap (\pi
  + \bm{z}/2).
\]
The singularity of the integrand is now located at $\bm{z}=0$.

\begin{figure}[hbt]
  \centering
  \includegraphics{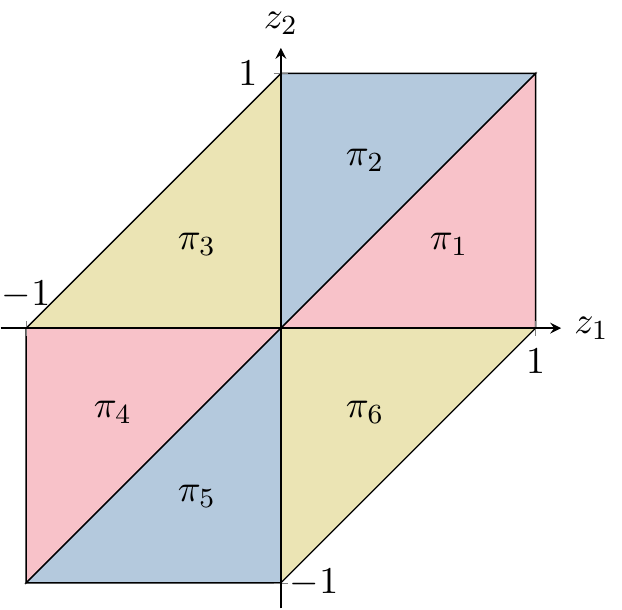}
  \caption{The domain of integration $\Pi$ is split into
    six triangles $\pi_i$.}\label{fig:domain_id}
\end{figure}

As shown in Figure~\ref{fig:domain_id}, we decompose $\Pi$ into six triangles
$\pi_i = \bm{A}_i \pi$ with
\[
  \bm{A}_1 = \begin{pmatrix} 1 & 0 \\ 0 & 1\end{pmatrix},\quad
  \bm{A}_2 = \begin{pmatrix} 0 & 1 \\ 1 & 0\end{pmatrix},\quad
  \bm{A}_3 = \begin{pmatrix} 0 & 1 \\ 1 & -1\end{pmatrix}
\]
and $\bm{A}_{i+3} = -\bm{A}_i$ for $i=1,2,3$. Thus, we obtain
\[
  I = \sum_{i=1}^6 \int\limits_{\pi} \int\limits_{\pi_{\bm{A}_i \bm{z}}}
  q(\bm{Z} + \bm{A}_i \bm{z}/2,
  \bm{Z} - \bm{A}_i \bm{z}/2) \, d\bm{Z}\, d\bm{z}.
\]
In the last step, we parametrise $\pi$ by ${(0,1)}^2$ using the Duffy
transformation 
\[
  \bm{z}(\bm{\eta}) = 
  \begin{pmatrix}
    \eta_1\\ \eta_1 \eta_2
  \end{pmatrix},
  \quad
  \bm{A}_i(\bm{\eta}) = \bm{A}_i \bm{z}(\bm{\eta})
\]
with Jacobian $\eta_1$ and conclude that $I$ has the representation 
\[
  I = \int\limits_{{(0,1)}^2} \eta_1 \sum_{i=1}^6\ \int\limits_{\pi_{\bm{A}_i
      (\bm{\eta})}} q(\bm{Z} + \bm{A}_i(\bm{\eta})/2, \bm{Z} -
  \bm{A}_i(\bm{\eta})/2)
  \, d\bm{Z}\, d\bm{\eta}.
\]
In Section~\ref{sec:calc}, we will see that the integrand is smooth since the
Jacobian $\eta_1$ cancels out the singularity of $q$. 

\subsection{Common edge}\label{subsec:reg_edge}
If the two triangles intersect at exactly one edge, we can proceed similarly
to the first case. Let $\chi_\sigma$ and $\chi_\tau$ be chosen in such a way that
the common edge is parametrised by
\[
  \chi_\tau (x_1,0) = \chi_\sigma(x_1,0), \quad x_1 \in (0,1).
\]
The kernel function is singular at this edge, i.e.\ at $(\bm{x},\bm{y})$ with
$x_2  = y_2 = 0$ and $x_1 = y_1$. The regularisation is
carried out via the mappings
\[
  \begin{gathered}
    \bm{A}_1 (\bm{\eta}) = \eta_1
    \begin{pmatrix}
      1\\ \eta_2\eta_4\\ 1-\eta_2\eta_3\\ \eta_2(1-\eta_3)
    \end{pmatrix},
    \ \,
    \bm{A}_2 (\bm{\eta}) = \eta_1
    \begin{pmatrix}
      1\\ \eta_2\\ 1-\eta_2\eta_3\eta_4\\ \eta_2\eta_3(1-\eta_4)
    \end{pmatrix},
    \ \,
    \bm{A}_3 (\bm{\eta}) = \eta_1
    \begin{pmatrix}
      1-\eta_2\eta_3\\ \eta_2(1-\eta_3)\\ 1\\ \eta_2\eta_3\eta_4
    \end{pmatrix},
    \\[1ex]
    \bm{A}_4 (\bm{\eta}) = \eta_1
    \begin{pmatrix}
      1-\eta_2\eta_3\eta_4\\ \eta_2\eta_3(1-\eta_4)\\ 1\\ \eta_2
    \end{pmatrix},
    \quad
    \bm{A}_5 (\bm{\eta}) = \eta_1
    \begin{pmatrix}
      1-\eta_2\eta_3\eta_4\\ \eta_2(1-\eta_3\eta_4)\\ 1\\ \eta_2\eta_3
    \end{pmatrix},
  \end{gathered}
\]
and reads
\[
  I = \int\limits_{{(0,1)}^4} \eta_1^3\eta_2^2
  \left(q(\bm{A}_1(\bm{\eta}))+\eta_3\sum_{i=2}^5 q(\bm{A}_i(\bm{\eta}))\right)
  d\bm{\eta}.
\]

\subsection{Common vertex}\label{subsec:reg_vert}
Let the origin in the reference domain be mapped to the common vertex, i.e.
\[
  \chi_\tau(0,0) = \chi_\sigma(0,0).
\]
By virtue of the mappings
\[
  \bm{A}_1 (\bm{\eta}) = \eta_1
  \begin{pmatrix}
    1\\ \eta_2\\ \eta_3\\\eta_3\eta_4
  \end{pmatrix},
  \quad
  \bm{A}_2 (\bm{\eta}) = \eta_1
  \begin{pmatrix}
    \eta_3\\ \eta_3\eta_4\\ 1\\ \eta_2
  \end{pmatrix},
\]
we obtain
\[
  I = \int\limits_{{(0,1)}^4} \eta_1^3 \eta_3 \left(q(\bm{A}_1(\bm{\eta})) +
    q(\bm{A}_2(\bm{\eta}))\right) 
  d\bm{\eta}.
\]

In summary, the regularisation yields integral representations with smooth
integrands on the unit cube. In this form, the integral can be approximated
efficiently by quadrature rules and the quadrature error decays exponentially
with the quadrature order.
\section{Calculation of integrals}\label{sec:calc}
Instead of applying quadrature rules directly, we calculate parts of the
regularised integrals analytically.

\subsection{Single layer potential}\label{subsec:slp}
With the discretisation provided in Section~\ref{sec:prelim}, the entries of the
single layer potential $\bm{V}$ are of the form
\begin{equation}\label{eq:slp}
  I = \int\limits_{\sigma} \int\limits_{\tau} \frac1{4\pi
    \seminorm{\bm{y}-\bm{x}}}\,dS(\bm{y})\, dS(\bm{x})
  = \frac{g_\tau g_\sigma}{4\pi} \int\limits_{\pi \times \pi}
  \frac1{\seminorm{\chi_\tau(\bm{y})-\chi_\sigma(\bm{x})}}
  \,d\bm{y}\,d\bm{x}.
\end{equation}
We proceed like in Section~\ref{sec:reg} and begin with the case of identical
elements.

\subsubsection{Identical Elements}\label{subsubsec:slp_id}
\begin{figure}[hbt]
  \centering
  \includegraphics{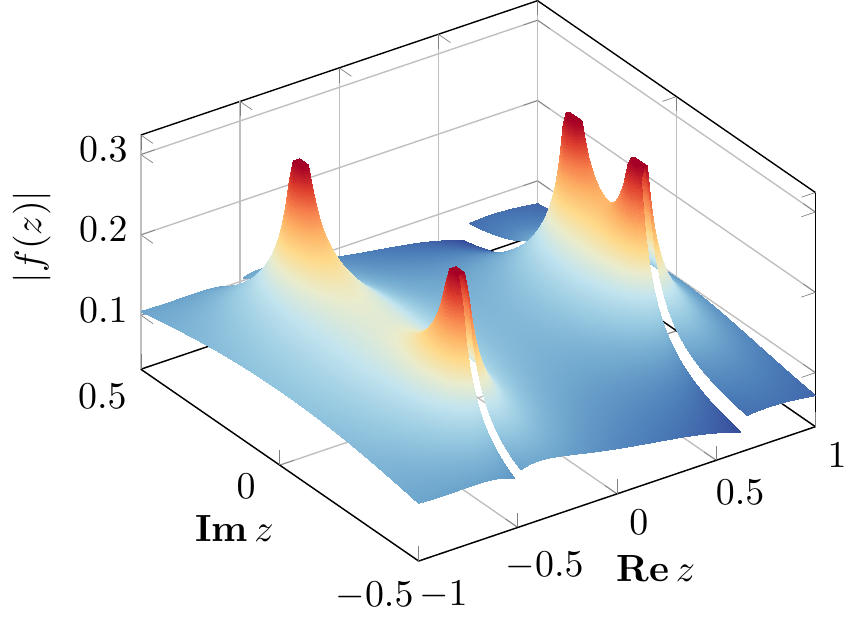}
  \caption{Visualisation of the integrand $f(z)$ of~\eqref{eq:plotid} in the
    complex plane.}\label{fig:slp_id_2D}
\end{figure}

Let $\bm{v}$ and $\bm{w}$ be the edges of $\tau=\sigma$ with starting point
$\bm{p}$. Then, the triangle is parametrised by
\[
  \chi_\sigma(\bm{y}) = \chi_\tau(\bm{y}) = \bm{p} + y_1\bm{v} + y_2 \bm{w}
\]
and the regularisation of~\eqref{eq:slp} reads
\[
  \begin{aligned}
    I &= \int\limits_{{(0,1)}^2} \eta_1 \sum_{i=1}^6\ \int\limits_{\pi_{\bm{A}_i
        (\bm{\eta})}} q(\bm{Z} + \bm{A}_i(\bm{\eta})/2, \bm{Z} -
    \bm{A}_i (\bm{\eta})/2) \, d\bm{Z}\, d\bm{\eta}\\[1ex]
    &= \frac{g_\tau^2}{2\pi} \int\limits_{{(0,1)}^2} \eta_1 \left(
      \frac{\seminorm{\pi_{\bm{A}_1(\bm{\eta})}}}
      {\seminorm{\eta_1 \eta_2 \bm{v} + \eta_1 \bm{w}}} +
      \frac{\seminorm{\pi_{\bm{A}_2 (\bm{\eta})}}}
      {\seminorm{\eta_1 \eta_2 \bm{w} + \eta_1 \bm{v}}} + 
      \frac{\seminorm{\pi_{\bm{A}_3 (\bm{\eta})}}}
      {\seminorm{\eta_1\eta_2 (\bm{w} + \bm{v}) - \eta_1\bm{w}}}
    \right)
    d\bm{\eta},
  \end{aligned}
\]
where the area is $\seminorm{\pi_{\bm{A}_i (\bm{\eta})}} = {(1-\eta_1)}^2/2$ for
$i=1,2,3$. Hence, we obtain
\begin{equation}\label{eq:plotid}
  I = \frac{g_\tau^2}{12\pi} \int\limits_{{(0,1)}} \left(
    \frac1{\seminorm{\eta_2 \bm{v} + \bm{w}}} +
    \frac1{\seminorm{\eta_2 \bm{w} + \bm{v}}} + 
    \frac1{\seminorm{\eta_2 (\bm{w} + \bm{v}) - \bm{w}}}
  \right) d\bm{\eta}.
\end{equation}
Figure~\ref{fig:slp_id_2D} depicts the complex continuation of the integrand for
concrete values of $\bm{v}$ and $\bm{w}$. It is smooth on the real axis, since
the edges are linearly independent, but has poles and branch cuts in the complex
domain. 

The three terms in the integrand are of the form
\[
  \frac1{\sqrt{\gamma + \beta \eta_2 + \alpha \eta_2^2 }},\quad
  \textup{with } \alpha>0,\ 4\alpha \gamma - \beta^2>0. 
\]
The anti-derivative is given by
\begin{equation}\label{eq:antid}
  F(\eta_2) = \frac1{\sqrt{\alpha}}\, \ln \left(2 \sqrt{\alpha} \sqrt{\gamma
      + \beta \eta_2 + \alpha \eta_2^2 } + 2 \alpha \eta_2 +
    \beta\right),
\end{equation}
see~\cite[Section~1.2.52.8]{Prudnikov1981}
and~\cite[Section~2.261]{Gradstein2015}. Thus, the integral reduces to 
\[
  I = \frac{g_\tau^2}{12\pi} {\Big[F_1(\eta_2)  + F_2(\eta_2) +
    F_3(\eta_2)\Big]}_0^1,
\]
where $F_i$ is $F$ with the parameters
\[
  \begin{aligned}
    \alpha_1 &= \seminorm{\bm{v}}^2,& \beta_1 &= 2 \bm{v} \cdot \bm{w},&
    \gamma_1 &= \seminorm{\bm{w}}^2,\\
    \alpha_2 &= \seminorm{\bm{w}}^2,& \beta_2 &= 2 \bm{w} \cdot \bm{v},&
    \gamma_2 &= \seminorm{\bm{v}}^2,\\
    \alpha_3 &= \seminorm{\bm{w}+\bm{v}}^2,& \beta_3 &= -2 (\bm{w}+\bm{v}) \cdot
    \bm{w},& \gamma_3 &= \seminorm{\bm{w}}^2.
  \end{aligned}
\]

\subsubsection{Common edge}\label{subsubsec:slp_edge}
\begin{figure}[hbt]
  \centering
  \includegraphics{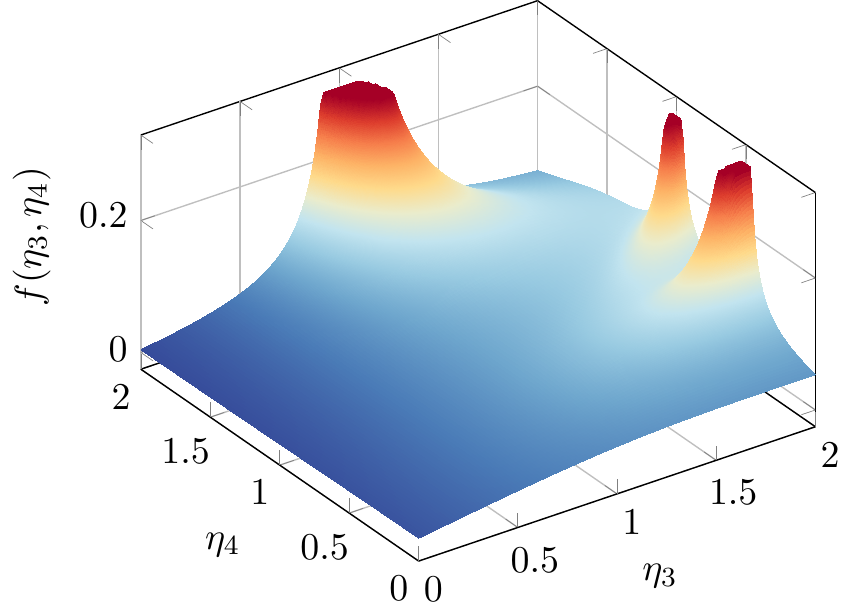}
  \caption{Visualisation of the integrand $f(\eta_3,\eta_4)$
    of~\eqref{eq:plot_slpedge}.}\label{fig:slp_edge_2D} 
\end{figure}

Let the reference mappings be given by
\[
  \chi_\tau(\bm{y}) = \bm{p} + y_1 \bm{v} + y_2 \bm{u},\quad
  \chi_\sigma(\bm{x}) = \bm{p} + x_1 \bm{v} + x_2 \bm{w},
\]
such that $\bm{v}$ is the common edge of $\sigma$ and $\tau$ starting from
$\bm{p}$. Then, the integral~\eqref{eq:slp} reduces to
\begin{equation}\label{eq:plot_slpedge}
  \begin{aligned}
    I &= \int\limits_{{(0,1)}^4} \eta_1^3\eta_2^2
    \left(q(\bm{A}_1(\eta))+\eta_3\sum_{i=2}^5 q(\bm{A}_i(\bm{\eta}))\right)
    d\bm{\eta}\\[1ex]
    &=  \frac{g_\tau g_\sigma}{24\pi} \int\limits_{{(0,1)}^2} \left(
      \frac1{\seminorm{\eta_3(\bm{u} + \bm{v}) + \eta_4 \bm{w} - \bm{u}}} +
      \frac{\eta_3}{\seminorm{\eta_4\eta_3(\bm{u}+\bm{v}) - \eta_3\bm{u} +
          \bm{w}}} \right.
    \\[1ex]
    &\quad +
    \frac{\eta_3}{\seminorm{\eta_4\eta_3(\bm{w}+\bm{v}) - \eta_3 \bm{w}
        + \bm{u}}} +
    \frac{\eta_3}{\seminorm{\eta_4\eta_3\bm{u} + \eta_3(\bm{w}+\bm{v}) -
        \bm{w}}}
    \\[1ex]
    &\quad+\left.
      \frac{\eta_3}{\seminorm{\eta_4\eta_3(\bm{w}+\bm{v}) + \eta_3 \bm{u} - \bm{w}}}
    \right) d \eta_3 d\eta_4
    \\[1ex]
    &= \frac{g_\tau g_\sigma}{24\pi} \sum_{i=1}^5 I_i.
  \end{aligned}
\end{equation}
In comparison to the previous case, the integrand is not necessarily smooth in
the real domain. When the two triangles lie in the same plane, it has poles as
seen in Figure~\ref{fig:slp_edge_2D}. However, they only occur outside of
${(0,1)}^2$ since the triangles do not overlap.

\subsubsection*{First integral $I_1$}
Let us introduce the variables
\[
  \bm{a} = \bm{w},\quad \bm{b} = \bm{v},\quad \bm{c} = \bm{u}+\bm{v}.
\]
We integrate with respect to $\eta_3$ by using~\eqref{eq:antid} and obtain for
the first integral
\[
  I_1 = \frac1{\seminorm{\bm{c}}} \int\limits_0^1 \ln \left(\frac{\seminorm{\eta_4
        \bm{a}+\bm{b}} \seminorm{\bm{c}} + (\eta_4 \bm{a}+\bm{b})\cdot
      \bm{c}}{\seminorm{\eta_4 \bm{a}+\bm{b}-\bm{c}} \seminorm{\bm{c}} +(\eta_4
      \bm{a}+\bm{b}-\bm{\bm{c}})\cdot \bm{c}}\right) d\eta_4. 
\]
Integration by parts leads to
\begin{equation}\label{eq:intedge}
  I_1 = \frac1{\seminorm{\bm{c}}} \ln\left(\frac{\seminorm{\bm{a}+\bm{b}}
      \seminorm{\bm{c}}  + (\bm{a}+\bm{b})\cdot \bm{c}}{\seminorm{\bm{a}+\bm{b}-
        \bm{c}}\seminorm{\bm{c}} + (\bm{a}+\bm{b}-\bm{c})\cdot \bm{c}}\right)
  - \frac1{\seminorm{\bm{c}}}
  \int\limits_0^1 \left( h_1(\eta_4) - h_0(\eta_4) \right) d\eta_4,
\end{equation}
where
\[
  \begin{aligned}
    h_0(\eta_4) &= \frac{(\eta_4 \bm{a}+\bm{b})\cdot \bm{b} + \seminorm{\eta_4
        \bm{a}+\bm{b}} \bm{b}\cdot \hat{\bm{c}}}{\seminorm{\eta_4
        \bm{a}+\bm{b}}^2 + \seminorm{\eta_4 \bm{a}+\bm{b}} (\eta_4
      \bm{a}+\bm{b}) \cdot \hat{\bm{c}}},\\[1ex]
    h_1(\eta_4) &= \frac{(\eta_4 \bm{a}+\bm{b}-\bm{c})\cdot (\bm{b}-\bm{c}) +
      \seminorm{\eta_4 \bm{a}+\bm{b}-\bm{c}} (\bm{b}-\bm{c}) \cdot
      \hat{\bm{c}}}{\seminorm{\eta_4 \bm{a}+\bm{b}-\bm{c}}^2 + \seminorm{\eta_4
        \bm{a}+\bm{b}-\bm{c}} (\eta_4 \bm{a}+\bm{b}-\bm{c}) \cdot \hat{\bm{c}}}
  \end{aligned}
\]
with $\hat{\bm{c}} = \bm{c} / \seminorm{\bm{c}}$. We note that $h_1$ coincides
with $h_0$ when $\bm{b}$ is replaced by $\bm{b}-\bm{c}$ and proceed with
integrating $h=h_0$. We follow the approach of~\cite[Appendix
C.2]{RjasanowSteinbach2007} and define 
\[
  p = \frac{\bm{a}\cdot \bm{b}}{\seminorm{\bm{a}}^2},\quad
  q^2 = \frac{\seminorm{\bm{b}}^2}{\seminorm{\bm{a}}^2} - p^2 \ge0
\]
such that
\[
  \seminorm{\eta_4 \bm{a}+\bm{b}} = \seminorm{\bm{a}} \sqrt{{(\eta_4+p)}^2+q^2}.
\]
If $q=0$ then $\bm{a}=p\, \bm{b}$ and the integral simplifies to
\[
  \int\limits h(\eta_4)\, d\eta_4 = p \ln (1 + 1/p).
\]
Otherwise, we have $q>0$ and the substitution
\[
  \eta_4 = -p + q \sinh (s), \quad d\eta_4 =q \cosh(s)ds,
\]
yields for the indefinite integral
\[
  \begin{aligned}
    \int\limits h(\eta_4)\,d\eta_4 &= \int\limits q \cosh(s) \frac{(-p+q\sinh(s))
      p \seminorm{\bm{a}}^2 + (q^2 +p^2) \seminorm{\bm{a}}^2 + \seminorm{\bm{a}}
      q \cosh(s) \bm{b}\cdot \hat{\bm{c}}}{q^2\seminorm{\bm{a}}^2 \cosh^2(s) + q
      \seminorm{\bm{a}} \cosh(s) [(-p+q\sinh(s))\bm{a}+\bm{b}]\cdot
      \hat{\bm{c}}} \,ds\\[1ex]
    &= q \int\limits \frac{p \seminorm{\bm{a}} \sinh(s) + q \seminorm{\bm{a}} +
      \cosh(s) \bm{b}\cdot \hat{\bm{c}}}{q \seminorm{\bm{a}}\cosh(s) + q\sinh(s)
      \bm{a}\cdot \hat{\bm{c}} +(-p\bm{a}+\bm{b})\cdot \hat{\bm{c}}}
    \,ds.
  \end{aligned}
\]
We use a variant of the Weierstra{\ss} substitution,
\[
  \tanh(s/2) = t,\quad \sinh(s) = \frac{2t}{1-t^2}, \quad \cosh(s) =
  \frac{1+t^2}{1-t^2},\quad ds = \frac2{1-t^2}\,dt,
\]
and obtain
\[
  \begin{aligned}
    \int\limits h(\eta_4)\,d\eta_4 &= \int\limits \frac{2q}{1-t^2}\ \frac{2p
      \seminorm{\bm{a}} t + q \seminorm{\bm{a}} (1-t^2) + (\bm{b}\cdot
      \hat{\bm{c}})\,(1+t^2)}{q\seminorm{\bm{a}}(1+t^2) + 2 q (\bm{a}\cdot 
      \hat{\bm{c}})\, t + ((-p\bm{a}+\bm{b})\cdot \hat{\bm{c}})\,(1-t^2)} \,dt\\[1ex]
    &= 2q \int\limits \frac{1}{1-t^2}\ \frac{\bm{b}\cdot \hat{\bm{c}} + q
      \seminorm{\bm{a}} + 2p \seminorm{\bm{a}} t + ((\bm{b}\cdot \hat{\bm{c}}) -
      q \seminorm{\bm{a}})\, t^2}{q\seminorm{\bm{a}} + (-p\bm{a}+\bm{b})\cdot
      \hat{\bm{c}} + 2q(\bm{a}\cdot \hat{\bm{c}})\, t + (q\seminorm{\bm{a}}
      + (p\bm{a}-\bm{b})\cdot \hat{\bm{c}})\,t^2} \,dt.
  \end{aligned}
\]
The integrand is now a rational function and we abbreviate it by
\[
  \frac{1}{1-t^2} \frac{\beta_0 + \beta_1 t + \beta_2 t^2}
  {\alpha_0 + \alpha_1 t + \alpha_2 t^2}.
\]
We decompose it into partial fractions,
\[
  \frac{\gamma_1}{1-t} + \frac{\gamma_2}{1+t} + \frac{\gamma_3 + \gamma_4 t}
  {\alpha_0 + \alpha_1 t + \alpha_2 t^2},
\]
where
\[
  \begin{aligned}
    \gamma_1 &= \frac12 \frac{\beta_0 + \beta_1 + \beta_2}
    {\alpha_0 + \alpha_1 + \alpha_2},&
    \gamma_2 &= \frac12 \frac{\beta_0 - \beta_1 + \beta_2}
    {\alpha_0 - \alpha_1 + \alpha_2},\\[1ex]
    \gamma_3 &= \beta_0 - (\gamma_1+\gamma_2)\alpha_0 ,&
    \gamma_4 &= \alpha_2(\gamma_1-\gamma_2). 
  \end{aligned}
\]
The first two terms yield
\[
  F(t) = \int\limits \left(\frac{\gamma_1}{1-t} + \frac{\gamma_2}{1+t}\right) dt
  = \gamma_2 \ln\seminorm{1+t} - \gamma_1 \ln\seminorm{1-t}.
\]
The third term depends on the discriminant $D =4\alpha_0\alpha_2-\alpha_1^2$ of
the denominator, which is non-negative due to
\[
  D=\seminorm{\det(\bm{a}|\bm{b}|\hat{\bm{c}})}^2 / \seminorm{\bm{a}}^2.
\]
For $D>0$ we have 
\[
  \begin{aligned}
    G(t) = \int\limits \frac{\gamma_3 + \gamma_4 t}
    {\alpha_0 + \alpha_1 t + \alpha_2t^2}\,dt &= \frac{\gamma_4}{2a_2} \ln
    \seminorm{\alpha_0 + \alpha_1 t + \alpha_2t^2}\\ 
    &\quad + \frac{2 \gamma_3 \alpha_2 - \gamma_4 \alpha_1}{\alpha_2 \sqrt{D}}
    \arctan \left(\frac{\alpha_1 + 2\alpha_2 t}{\sqrt{D}} \right)
  \end{aligned}
\]
and for $D=0$
\[
  G(t) = \int\limits \frac{\gamma_3 + \gamma_4 t}
  {\alpha_0 + \alpha_1 t + \alpha_2t^2}\,dt = \frac{\gamma_4}{\alpha_2}
  \ln \seminorm{t + \frac{\alpha_1}{2\alpha_2}}
  - \frac{2 \gamma_3 \alpha_2 - \gamma_4 \alpha_1}
  {\alpha_2 (2\alpha_2t + \alpha_1)}. 
\]
We resubstitute
\[
  t = \tanh \left[\frac12
    \arcsinh \left(\frac{p+\eta_4}{q}\right)\right]
  = \frac{\sinh\left[\arcsinh\left(\frac{p+\eta_4}{q}\right)\right]}{1 +
    \cosh\left[\arcsinh\left(\frac{p+\eta_4}{q}\right)\right]} 
  = \frac{p + \eta_4}{q + \sqrt{{(p+\eta_4)}^2+q^2}}, \quad
\]
and set
\[
  t_0 = \frac{p}{q + \sqrt{p^2+q^2}}, \quad
  t_1 = \frac{p + 1}{q + \sqrt{{(p+1)}^2+q^2}}.
\]
Finally, we obtain
\begin{equation}\label{eq:antid2}
  \int\limits_0^1 h(\eta_4)\,d\eta_4 = 2q \left( F(t_1) - F(t_0) + G(t_1) - G(t_0)
  \right).
\end{equation}
Note that the value of the integral only depends on the vectors $\bm{a}, \bm{b},
\bm{c}$. Since it is of importance for the other cases as well, we abbreviate it
by
\[
  H(\bm{a},\bm{b},\bm{c}) = \int\limits_0^1 h(\eta_4)\,d\eta_4.
\]
We conclude that $I_1$ can be expressed in closed form as
\begin{equation}\label{eq:edge1}
  I_1 = \frac1{\seminorm{\bm{c}}}
  \ln\left(\frac{\seminorm{\bm{a}+\bm{b}} \seminorm{\bm{c}} +
      (\bm{a}+\bm{b})\cdot \bm{c}}{\seminorm{\bm{a}+\bm{b}-\bm{\bm{c}}}
      \seminorm{\bm{c}}+(\bm{a}+\bm{b}-\bm{c})\cdot \bm{c}}\right)
  - \frac1{\seminorm{\bm{c}}} \left(H(\bm{a},\bm{b}-\bm{c},\bm{c}) -
    H(\bm{a},\bm{b},\bm{c}) \right)
\end{equation}
with $\bm{a} = \bm{w}$, $\bm{b} = \bm{v}$, $\bm{c} = \bm{u}+\bm{v}$.

\subsubsection*{Remaining integrals}
For the remaining integrals, we integrate with respect to the fourth variable
firstly. With
\[
  \bm{a} = \bm{v},\quad \bm{b} = \bm{w},\quad \bm{c} = \bm{u}+\bm{v},
\]
we have
\[
  I_2 = \frac1{\seminorm{\bm{c}}} \int\limits_0^1 \ln \left(\frac{\seminorm{\eta_3
        \bm{a}+\bm{b}} \seminorm{\bm{c}} + (\eta_3 \bm{a}+\bm{b})\cdot
      \bm{c}}{\seminorm{\eta_3 (\bm{a}-\bm{c})+\bm{b}} \seminorm{\bm{c}}
      +(\eta_3 (\bm{a}-\bm{c})+\bm{b})\cdot \bm{c}}\right) d\eta_3.
\]
Integration by parts yields an expression almost identical to~\eqref{eq:intedge},
\[
  I_2 = \frac1{\seminorm{\bm{c}}} \ln\left(\frac{\seminorm{\bm{a}+\bm{b}}
      \seminorm{\bm{c}} + (\bm{a}+\bm{b})\cdot \bm{c}}{\seminorm{\bm{a}-\bm{c}
        +\bm{b}} \seminorm{\bm{c}} + (\bm{a}-\bm{c}+\bm{b})\cdot \bm{c}}\right)
  - \frac1{\seminorm{\bm{c}}}
  \int\limits_0^1 \left( h_1(\eta_3) - h_0(\eta_3) \right) d\eta_3,
\]
where $h_0$ and $h_1$ are given by
\[
  \begin{aligned}    
    h_0(\eta_3) &= \frac{(\eta_3 \bm{a}+\bm{b})\cdot \bm{b} + \seminorm{\eta_3
        \bm{a}+\bm{b}} \bm{b}\cdot \hat{\bm{c}}}{\seminorm{\eta_3
        \bm{a}+\bm{b}}^2 + \seminorm{\eta_3 \bm{a}+\bm{b}} (\eta_3
      \bm{a}+\bm{b}) \cdot \hat{\bm{c}}},\\[1ex]
    h_1(\eta_3) &= \frac{(\eta_3 (\bm{a}-\bm{c})+\bm{b})\cdot \bm{b} +
      \seminorm{\eta_3 (\bm{a}-\bm{c})+\bm{b}} \bm{b}\cdot
      \hat{\bm{c}}}{\seminorm{\eta_3 (\bm{a}-\bm{c})+\bm{b}}^2 +
      \seminorm{\eta_3 (\bm{a}-\bm{c})+\bm{b}} (\eta_3 (\bm{a}-\bm{c})+\bm{b})
      \cdot \hat{\bm{c}}}. 
  \end{aligned}
\]
Thus, $I_2$ can be computed analogously to~\eqref{eq:edge1} by
\begin{equation}\label{eq:edgej}
  \frac1{\seminorm{\bm{c}}} \ln\left(\frac{\seminorm{\bm{a}+\bm{b}}
      \seminorm{\bm{c}} + (\bm{a}+\bm{b})\cdot \bm{c}}{\seminorm{\bm{a}-\bm{c}
        +\bm{b}} \seminorm{\bm{c}} + (\bm{a}-\bm{c}+\bm{b})\cdot
      \bm{c}}\right)
  - \frac1{\seminorm{\bm{c}}} \left(H(\bm{a}-\bm{c},\bm{b},\bm{c}) -
    H(\bm{a},\bm{b},\bm{c}) \right).
\end{equation}
Because this applies to the other integrals as well, we only list the parameters
$\bm{a}$, $\bm{b}$ and $\bm{c}$ in Table~\ref{table:slp_paredge}.

\begin{table}[htb]
  \centering
  \begin{tabular}{c|c|c|c}
    $j$ & $\bm{a}$ & $\bm{b}$ & $\bm{c}$\\
    \hline\hline
    $2$ & $\bm{v}$ & $\bm{w}$ & $\bm{u}+\bm{v}$\\
    \hline
    $3$ & $\bm{v}$ & $\bm{u}$ & $\bm{w}+\bm{v}$\\
    \hline
    $4$ & $\bm{u}+\bm{v}+\bm{w}$ & $-\bm{w}$ & $\bm{u}$\\
    \hline
    $5$ & $\bm{u}+\bm{v}+\bm{w}$ & $-\bm{w}$ & $\bm{v}+\bm{w}$\\
    \hline 
  \end{tabular}
  \caption{Values for $\bm{a},\bm{b},\bm{c}$ in~\eqref{eq:edgej} to compute
    $I_j$.}\label{table:slp_paredge} 
\end{table}

\subsubsection{Common vertex}\label{subsubsec:slp_vert}
We consider the configuration
\[
  \chi_\tau(\bm{y}) = \bm{p} + y_1 \bm{u}_1 + y_2 \bm{u}_2,\quad
  \chi_\sigma(\bm{x}) = \bm{p} + x_1 \bm{v}_1 + x_2 \bm{v}_2
\]
with common vertex $\bm{p}$. Then, integration with respect to $\eta_1$ results
in
\[
  \begin{aligned}
    I &=\frac{g_\tau g_\sigma}{12\pi} \! \int\limits_{{(0,4)}^3} \! \left(
      \frac{\eta_3}{\seminorm{\eta_3 \bm{u}_1 + \eta_3 \eta_4 \bm{u}_2 - \bm{v}_1-
          \eta_2 \bm{v}_2}} + \frac{\eta_3}{\seminorm{\bm{u}_1 + \eta_2 \bm{u}_2 -
          \eta_3 \bm{v}_1 - \eta_3\eta_4 \bm{v}_2}} \right)
    d\eta_2 d\eta_3 d\eta_4\\[1ex]
    &= \frac{g_\tau g_\sigma}{12\pi} \left(I_1+I_2\right).
  \end{aligned}
\]
We only consider $I_1$, since $I_2$ is obtained by swapping $\bm{u}_i$ and
$\bm{v}_i$. Similar to the previous section, we introduce the variables 
\[
  \bm{a}=\bm{u}_1+\bm{u}_2,\quad \bm{b}(\eta_2)=-\bm{v}_1-\eta_2 \bm{v}_2,\quad
  \bm{c}=\bm{u}_2,
\]
and integrate with respect to $\eta_4$ to obtain
\[
  I_1= \frac1{\seminorm{\bm{c}}} \int\limits_{{(0,1)}^2} \ln
  \left(\frac{\seminorm{\eta_3 \bm{a}+\bm{b}(\eta_2)} \seminorm{\bm{c}} +
      (\eta_3 \bm{a}+\bm{b}(\eta_2))\cdot \bm{c}}{\seminorm{\eta_3
        (\bm{a}-\bm{c})+\bm{b}(\eta_2)} \seminorm{\bm{c}} +(\eta_3
      (\bm{a}-\bm{c})+\bm{b}(\eta_2))\cdot \bm{c}}\right) 
  d\eta_3\, d\eta_2.
\]
We insert Formula~\eqref{eq:edgej} for the inner integral, which yields
\[
  \begin{aligned}
    I_1&=\frac1{\seminorm{\bm{c}}} \int\limits_0^1\ln\left(
      \frac{\seminorm{\bm{a}+\bm{b}(\eta_2)}\seminorm{\bm{c}} +
        (\bm{a}+\bm{b}(\eta_2))\cdot \bm{c}}{\seminorm{\bm{a}-\bm{c} 
          +\bm{b}(\eta_2)} \seminorm{\bm{c}} +
        (\bm{a}-\bm{c}+\bm{b}(\eta_2))\cdot \bm{c}}\right) d\eta_2\\[1ex]
    &\quad-\frac1{\seminorm{\bm{c}}} \int\limits_0^1 \left(
      H(\bm{a}-\bm{c},\bm{b}(\eta_2),\bm{c}) -
      H(\bm{a},\bm{b}(\eta_2),\bm{c}) \right) d\eta_2. 
  \end{aligned}
\]
The first integral can be written in the form of $I_1$ from
Section~\ref{subsubsec:slp_edge}, i.e.
\[
  \frac1{\seminorm{\tilde{\bm{c}}}} \int\limits_0^1\ln\left(
    \frac{\seminorm{\eta_2\tilde{\bm{a}}+\tilde{\bm{b}}}
      \seminorm{\tilde{\bm{c}}} + (\eta_2\tilde{\bm{a}}+\tilde{\bm{b}})\cdot
      \tilde{\bm{c}}}
    {\seminorm{\eta_2\tilde{\bm{a}}+\tilde{\bm{b}}-\tilde{\bm{c}}}
      \seminorm{\tilde{\bm{c}}} + (\eta_2\tilde{\bm{a}}+\tilde{\bm{b}}-
      \tilde{\bm{c}})\cdot \tilde{\bm{c}}}\right) d\eta_2 
\]
with $\tilde{\bm{a}}=-\bm{v}_2$, $\tilde{\bm{b}}=\bm{u}_1+\bm{u}_2-\bm{v}_1$,
$\tilde{\bm{c}}= \bm{u}_2$, and its value is hence given by~\eqref{eq:edge1}.
Because it is not possible to integrate the remaining integral analytically, we
approximate it numerically with a quadrature rule
\[
  \begin{aligned}
    \int\limits_0^1 \left( H(\bm{a}-\bm{c},\bm{b}(\eta_2),
      \bm{c})-H(\bm{na},\bm{b}(\eta_2),\bm{c}) \right) d\eta_2& \approx \\
    &\hspace{-10em} \sum_{i=1}^n \omega_i \left( 
      H(\bm{a}-\bm{c},\bm{b}(\eta^{(i)}),\bm{c})-H(\bm{a},\bm{b}(\eta^{(i)}),\bm{c})
    \right)
  \end{aligned}
\]
with weights $\omega_i>0$ and nodes $\eta^{(i)}\in [0,1]$.

\subsubsection{Far-field}\label{subsubsec:slp_far}
Although the far-field does not constitute a singular case, the analytical
formulae are still applicable. Let the elements be given by
\[
  \chi_\tau(\bm{y}) = \bm{p}_1 + y_1 \bm{u}_1 + y_2 \bm{u}_2,\quad
  \chi_\sigma(\bm{x}) = \bm{p}_2 + x_1 \bm{v}_1 + x_2 \bm{v}_2,
\]
and set $\bm{p}=\bm{p}_1-\bm{p}_2$. Analogously to the previous cases, we pull
the region of integration back to ${(0,1)}^4$ by
\[
  \bm{A}: {(0,1)}^4 \to \pi \times \pi,\quad
  \bm{A}(\bm{\eta}) = \begin{pmatrix}
    \eta_1\\ \eta_1 \eta_2\\ \eta_3\\ \eta_3 \eta_4
  \end{pmatrix},
\]
leading to
\[
  I =\frac{g_\tau g_\sigma}{4\pi}\int\limits_{{(0,1)}^4}
  \frac{\eta_1\eta_3}{\seminorm{\bm{p}+\eta_3 \bm{u}_1 + \eta_3 \eta_4
      \bm{u}_2 - \eta_1 \bm{v}_1- \eta_1 \eta_2 \bm{v}_2}}\,
  d\bm{\eta}.
\]
Of the four iterated integrals, we compute two analytically and two by numerical
quadrature, e.g.
\[
  I \approx \sum_{k,\ell=1}^n \omega_k\, \omega_\ell \int\limits_{{(0,1)}^2}
  \frac{\eta^{(k)}\,\eta^{(\ell)}}{\seminorm{\bm{p}+\eta^{(\ell)} \bm{u}_1
      + \eta_4 \eta^{(\ell)} \bm{u}_2 - \eta^{(k)} \bm{v}_1-
      \eta_2 \eta^{(k)} \bm{v}_2}}\,
  d\eta_2\,d\eta_4,
\]
where the two-dimensional integral is calculated analytically
using~\eqref{eq:edge1} with 
\[
  \bm{a} = \eta^{(k)} \bm{v}_2,\quad
  \bm{b} = \bm{p} - \eta^{(k)} \bm{v}_1 + \eta^{(\ell)}(\bm{u}_1+\bm{u}_2),\quad
  \bm{c} = \eta^{(\ell)} \bm{u}_2.
\]

\subsection{Double layer potential}\label{subsec:dlp}

For the double layer potential $\bm{K}$, we need to compute integrals of the
form 
\begin{equation}\label{eq:dlp}
  J = \int\limits_{\sigma} \int\limits_{\tau}
  \frac{(\bm{x}-\bm{y})\cdot \bm{n}}{4\pi\seminorm{\bm{x}-\bm{y}}^3}
  \, \varphi(\bm{y})\,dS(\bm{y})\, dS(\bm{x})
  = \frac{g_\tau g_\sigma}{4\pi}\! \int\limits_{\pi \times \pi}
  \frac{\left(\chi_\sigma(\bm{x})-\chi_\tau(\bm{y})\right)
    \cdot\bm{n}}{\seminorm{\chi_\sigma(\bm{x})-\chi_\tau(\bm{y})}^3}
  \,\varphi(\chi_\tau(\bm{y})) \,d\bm{x}\,d\bm{y},
\end{equation}
where $\bm{n}$ is the outer unit normal vector at $\tau$ and $\varphi\in
S_h^1(\tau)$, i.e.
\[
  \varphi(\chi_\tau(\bm{y})) = a_0 + a_1 y_1 + a_2 y_2
\]
with coefficients $a_0,a_1,a_2\in\mathbb{R}$.

\subsubsection{Identical elements}\label{subsubsec:dlp_id}
For identical elements $\sigma=\tau$, we simply have $J = 0$ due to
\[
  (\bm{y}-\bm{x}) \cdot \bm{n} = 0,\quad \textup{for } \bm{x},\bm{y}\in\tau.
\]

\subsubsection{Common edge}\label{subsubsec:dlp_edge}
\begin{figure}[hbt]
  \centering
  \includegraphics{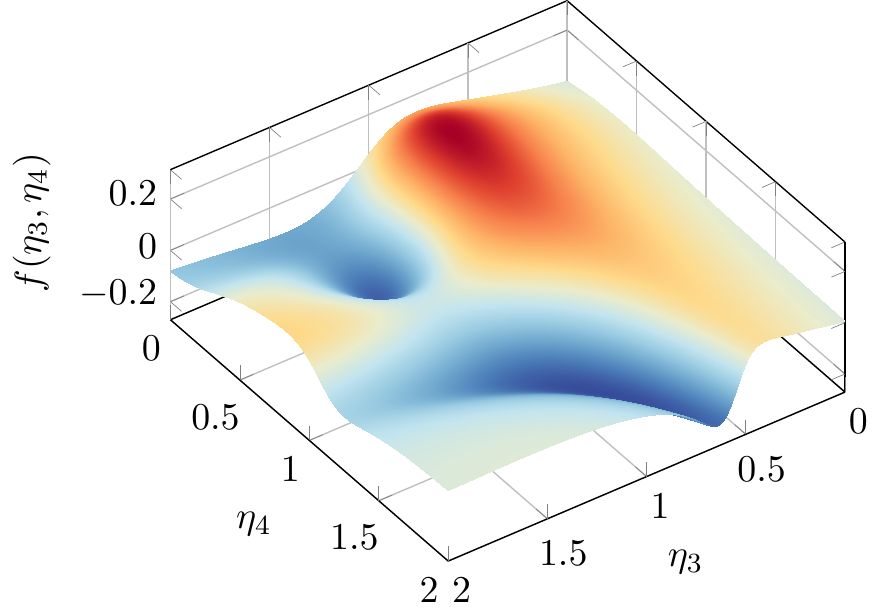}
  \caption{Visualisation of the integrand $f(\eta_3,\eta_4)$
    of~\eqref{eq:plot_dlpedge}.}\label{fig:dlp_edge_2D}
\end{figure}

We assume that the triangles are parametrised by
\[
  \chi_\tau(\bm{y}) = \bm{p} + y_1 \bm{v} + y_2 \bm{u},\quad
  \chi_\sigma(\bm{x}) = \bm{p} + x_1 \bm{v} + x_2 \bm{w}.
\]
Applying the regularisation to $J$ and integrating with respect to $\eta_1$ and
$\eta_2$ leads to 
\begin{equation}\label{eq:plot_dlpedge}
  J = \frac{g_\sigma g_\tau}{4\pi}\, \bm{w}\cdot\bm{n} \sum_{i=1}^5 J_i,
\end{equation}
where the integrals $J_i$ are given by
\[
  \begin{aligned}
    J_1 &= \int\limits_{{(0,1)}^2}
    \frac{\eta_4}
    {\seminorm{\eta_3(\bm{u}+\bm{v})+\eta_4\bm{w}-\bm{u}}^3}
    \left(c_0 +c_1\eta_3+c_2(1-\eta_3)\right) d\eta_3\,d\eta_4,\\[1ex]
    J_2 &= \int\limits_{{(0,1)}^2}
    \frac{\eta_3}
    {\seminorm{\eta_3\eta_4(\bm{u}+\bm{v})-\eta_3\bm{u}+\bm{w}}^3}
    \left(c_0+c_1\eta_3\eta_4+c_2\eta_3(1-\eta_4)\right) d\eta_3\,d\eta_4,\\[1ex]
    J_3 &= \int\limits_{{(0,1)}^2}
    \frac{\eta^2_3(1-\eta_4)}
    {\seminorm{\eta_3\eta_4(\bm{w}+\bm{v})-\eta_3\bm{w}+\bm{u}}^3}
    \left(c_0+c_2\right) d\eta_3\,d\eta_4,\\[1ex]
    J_4 &= \int\limits_{{(0,1)}^2}
    \frac{\eta_3(1-\eta_3)}
    {\seminorm{\eta_3\eta_4\bm{u}+\eta_3(\bm{w}+\bm{v})-\bm{w}}^3}
    \left(c_0+c_2\eta_3\eta_4\right) d\eta_3\,d\eta_4,\\[1ex]
    J_5 &= \int\limits_{{(0,1)}^2}
    \frac{\eta_3(1-\eta_3\eta_4)}
    {\seminorm{\eta_3\eta_4(\bm{w}+\bm{v})+\eta_3\bm{u}-\bm{w}}^3}
    \left(c_0+c_2\eta_3\right) d\eta_3\,d\eta_4.
  \end{aligned}
\]
with $c_0 = a_0 / 2 + a_1/3$, $c_1=-a_1/6$ and $c_2=a_2/6$. In contrast to the
respective case of the single layer potential, the integrand of $J$ is always
smooth in the real domain as shown in Figure~\ref{fig:dlp_edge_2D}. Indeed, if
the two triangles lie in the same plane, then $J=0$ due to $\bm{w}\cdot\bm{n}=0$.

In the following, we derive an analytic expression for the integrals by the
example of $J_1$. We define
\[
  R(\eta_3,\eta_4) = \sqrt{\gamma(\eta_3)+\beta(\eta_3) \eta_4 +\alpha
    \eta_4^2}
\]
with
\[
  \gamma(\eta_3) = \seminorm{\eta_3(\bm{u}+\bm{v})-\bm{u}}^2,\quad
  \beta(\eta_3)=2(\eta_3(\bm{u}+\bm{v})-\bm{u})\cdot\bm{w},\quad
  \alpha = \seminorm{\bm{w}}^2,
\]
such that $R^3$ equals the denominator. Denoting the discriminant by
$D=4\alpha\gamma - \beta^2$, we integrate with respect to $\eta_4$
using~\cite[Section~2.264]{Gradstein2015},  
\[
  \begin{aligned}
    J_1 &= -2\int\limits_0^1(c_0+c_2+(c_1-c_2)\eta_3)
    {\left[\frac{2\gamma(\eta_3)+\beta(\eta_3)\eta_4}
        {D(\eta_3)R(\eta_3,\eta_4)}\right]}_0^1 d\eta_3\\[1ex]
    &=-2\int\limits_0^1(c_0+c_2+(c_1-c_2)\eta_3)
    \left[\frac{2\gamma(\eta_3)+\beta(\eta_3)}{D(\eta_3)
        \sqrt{\gamma(\eta_3)+\beta(\eta_3)+\alpha}} -
      \frac{2\gamma(\eta_3)}{ D(\eta_3)\sqrt{\gamma(\eta_3)}}\right] d\eta_3.
  \end{aligned}
\]
Hence, the integral reduces to
\[
  J_1 = \int\limits_0^1 \left( h^{(1)}_1(\eta_3)\,d\eta_3 -  h^{(1)}_0(\eta_3)
  \right) \,d\eta_3
\]
with $h^{(1)}_0$ and $h^{(1)}_1$ of the form
\begin{equation}\label{eq:dlp_hform}
  h(\eta_3) = \frac{4P(\eta_3)}{D(\eta_3)\sqrt{Q(\eta_3)}}, 
\end{equation}
where $P(\eta_3)$ is a cubic and $Q(\eta_3)$ a positive quadratic polynomial
respectively,
\[
  P(\eta_3) = p_0+p_1\eta_3+p_2\eta_3^2+p_3\eta_3^3,\quad
  Q(\eta_3) = q (q_0+q_1\eta_3+\eta_3^2)>0.
\]
In order to find the anti-derivative of $h$, we write $D(\eta_3)$ as
\[
  D(\eta_3) = 4d(d_0 + d_1 \eta_3 + \eta_3^2)
\]
and decompose into partial fractions,
\[
  \begin{aligned}
    h(\eta_3) &= \frac{1}{d \sqrt{Q(\eta_3)}} \biggl(
    p_2-d_1p_3+p_3\eta_3 +\\
    &\hspace{6.5em}+ \frac{p_0-d_0p_2+d_0d_1p_3+(p_1-d_1p_2-d_0p_3+
      d_1^2p_3)\eta_3}{d_0+d_1\eta_3+\eta_3^2}
    \biggr).
  \end{aligned}
\]
We are familiar with the first term and recall that
\[
  \int\limits \frac{1}{\sqrt{q_0 + q_1
      \eta_3 + \eta_3^2}}\, d\eta_3 = F(\eta_3)
  = \ln \left(2 \sqrt{q_0 + q_1
      \eta_3 + \eta_3^2 } + 2 \eta_3 + q_1 \right) 
\]
and~\cite[Section~2.264]{Gradstein2015} also gives
\[
  \int\limits \frac{\eta_3}{\sqrt{q_0 + q_1 \eta_3 +
      \eta_3^2}}\,d\eta_3 = \sqrt{q_0 + q_1 \eta_3 +
    \eta_3^2} - \frac{q_1}{2} F(\eta_3).
\]
For the remaining term, we abbreviate the constants in the numerator by
\[
  n = p_0-d_0p_2+d_0d_1p_3,\quad
  m = p_1-d_1p_2-d_0p_3+d_1^2p_3.
\]
Following~\cite[Chapter~3]{FichtenholzII2003}, we substitute
\[
  \eta_3 = \frac{\nu+\mu t}{1+t},\quad
  d\eta_3 = \frac{\mu-\nu}{{(1+t)}^2}\,dt,
\]
where $\mu$ and $\nu$ are the real and distinct solutions of the quadratic
equation 
\[
  (d_1-q_1) z^2 + 2(d_0-q_0)z + (q_1 d_0 - d_1 q_0) = 0.
\]
In this way, the linear terms of the denominator vanish and we obtain
\[
  \omega \int\limits \sgn (1+t) \frac{n+m\nu + (n+m\mu)t}{{(\lambda+t^2)
      \sqrt{\varkappa + t^2}}}\,dt
\]
with
\[
  \begin{gathered}
    \omega = \frac{\mu-\nu}{(\mu^2+d_1\mu+d_0)\sqrt{\mu^2+q_1\mu+q_0}},\\[1ex]
    \lambda=\frac{\nu^2+d_1\nu+d_0}{\mu^2+d_1\mu+d_0},\qquad
    \varkappa=\frac{\nu^2+q_1\nu+q_0}{\mu^2+q_1\mu+q_0}.
  \end{gathered}
\]
We have
\[
  \int\limits \frac{(n+m\mu)t}{{(\lambda+t^2) \sqrt{\varkappa + t^2}}}\,dt =
  (n+m\mu)\int\limits \frac{1}{{\lambda-\varkappa+s^2}}\,ds  
\]
by means of $s_1 = \sqrt{\varkappa + t^2}$. For the remaining term, we use the
substitution 
\[
  s_0 = \frac{t}{\sqrt{\varkappa + t^2}},\quad \frac{dt}{\sqrt{\varkappa+t^2}} =
  \frac{ds_0}{1-s_0^2},
\]
such that
\[
  \lambda + t^2 = \frac{\lambda+(\varkappa - \lambda) s_0^2}{1-s_0^2}
\]
and
\[
  \int\limits \frac{n+m\nu}{{(\lambda+t^2) \sqrt{\varkappa + t^2}}}\,dt = (n+m\nu)
  \int\limits \frac{1}{\lambda+(\varkappa - \lambda) s_0^2}\,ds_0.
\]
The integrals are of the form
\[
  \int\limits \frac{1}{\rho+s^2}\,ds,\quad \textup{where }
  \rho = \lambda-\varkappa \textup{ or }
  \rho = \frac{\lambda}{\varkappa-\lambda},
\]
and the anti-derivative depends on the sign of $\rho$, 
\[
  \int\limits\frac{1}{{\rho+s^2}}\,ds = G(\rho,s) = \left\{
    \begin{aligned}
      &\frac{\arctan\left(s/\sqrt{\rho}\right)}
      {\sqrt{\rho}}, && \rho>0,\\
      &-\frac{1}{s}, && \rho=0,\\
      &\frac{1}{2\sqrt{-\rho}} \ln \seminorm{\frac{s-\sqrt{-\rho}}
        {s+\sqrt{-\rho}}}, && \rho<0,
    \end{aligned}
  \right.
\]
see~\cite[Section~2.103]{Gradstein2015}. The value of the indefinite integral is
therefore 
\begin{equation}\label{eq:antid3}
  \begin{aligned}
    \int\limits \frac{n+m\,\eta_3}{(d_0+d_1\eta_3+\eta_3^2)\sqrt{q_0 + q_1 \eta_3
        +\eta_3^2}}\,d\eta_3 &=\\[2ex]
    &\hspace{-15em}\omega \sgn (1+t) \left[
      \frac{n+m\nu}{\varkappa-\lambda}\,
      G\!\left(\frac{\lambda}{\varkappa-\lambda}, s_0(t)\right)+
      (n+m\mu)\,G\!\left(\lambda-\varkappa, s_1(t)\right)
    \right]
  \end{aligned}
\end{equation}
and we write $H(t)$ for short.
It remains to determine the domain of the integral in terms of $t$. We
resubstitute 
\[
  t^{-1}(0,1) = \left\{
    \begin{aligned}
      & (t_0,t_1), && t_0 > -1 \textup{ or } t_1 < -1. \\
      & (-\infty,t_0) \cup (t_1,\infty), && t_0<-1<t_1,
    \end{aligned}
  \right\}
\]
with $t_0 = \min\{t(0),t(1)\}$, $t_1 = \max\{t(0),t(1)\}$ and
\[
  t(0) = -\frac{\nu}{\mu}, \quad t(1) = -\frac{\nu-1}{\mu-1}.
\]
The second case requires the calculation of improper integrals. We have
\[
  \lim_{t \to \pm \infty} G(\rho,s_0(t)) = \lim_{s\to\pm 1} G(\rho,s)
  = G(\rho,\pm1)
\]
and
\[
  \lim_{t \to \pm \infty} G(\rho,s_1(t)) = \lim_{s\to\infty} G(\rho,s) = \left\{
    \begin{aligned}
      &\frac{\pi}{2 \sqrt{\rho}}, && \rho>0,\\
      &0, && \rho\le0.
    \end{aligned}
  \right.
\]
We combine the results and obtain
\begin{equation}\label{eq:dlp_intedge}
  \begin{aligned}
    \int\limits_0^1 h(\eta_3)\,d\eta_3 &= \frac{1}{d\sqrt{q}} {\left[
        p_3\sqrt{q_0+q_1\eta_3+\eta_3^2}+\left(p_2-d_1p_3-p_3\frac{q_1}{2}\right)
        F(\eta_3)\right]}_0^1\\[1.5ex]
    &\quad + \left\{
      \begin{aligned}
        & H(t_1) - H(t_0), && t_0 > -1 \textup{ or } t_1 < -1, \\
        & H(t_0) - \lim_{t\to-\infty} H(t) + \lim_{t\to\infty} H(t) - H(t_1), &&
        t_0<-1<t_1. 
      \end{aligned}
    \right.
  \end{aligned}
\end{equation}
Finally, the integral $J_1$ can now be computed by applying the formula to the
integrals of $h^{(1)}_1$ and $h^{(1)}_0$. The other $J_i$ can be calculated in
the same way. After integrating with respect to $\eta_4$, we have
\[
  J_i = \int\limits_0^1 \left( h^{(i)}_1(\eta_3)\,d\eta_3 - 
  h^{(i)}_0(\eta_3) \right) d\eta_3,
\]
where $h^{(i)}_0$ and $h^{(i)}_1$ are again of the form~\eqref{eq:dlp_hform}.
Thus, we conclude that all integrals are expressible analytically
via~\eqref{eq:dlp_intedge} in terms of certain parameters, which are listed in
Tables~\ref{table:dlp_paredge_q},~\ref{table:dlp_paredge_d}
and~\ref{table:dlp_paredge_p}.

\subsubsection{Common vertex}\label{subsubsec:dlp_vert}
Like in Section~\ref{subsubsec:slp_vert}, let $\bm{u}_1,\bm{u}_2$ be the
edges of $\tau$ and $\bm{v}_1,\bm{v}_2$ the edges of $\sigma$. We integrate
with respect to $\eta_1$ and obtain
\[
  \begin{aligned}
    J &= \frac{g_\tau g_\sigma}{4\pi} \left(J_1+J_2\right)\\
    &= \frac{g_\tau g_\sigma}{4\pi} \left(\ \int\limits_{{(0,4)}^3}
      \frac{\eta_3(\bm{v}_1+\eta_2\bm{v}_2)\cdot\bm{n}}
      {\seminorm{\eta_3\bm{u}_1+\eta_3\eta_4 \bm{u}_2 - \bm{v}_1-
          \eta_2 \bm{v}_2}^3} \left(c_0+c_1\eta_3+c_2\eta_3\eta_4 \right)
      d\eta_2\,d\eta_3\,d\eta_4 \right.\\
    &\hspace{5em} \left.+
      \int\limits_{{(0,4)}^3}\frac{\eta_3^2(\bm{v}_1+\eta_4\bm{v}_2)\cdot\bm{n}}
      {\seminorm{\bm{u}_1 + \eta_2 \bm{u}_2 -\eta_3 \bm{v}_1 - \eta_3\eta_4
          \bm{v}_2}^3} \left(c_0+c_1+c_2\eta_2\right) d\eta_2\,d\eta_3\,d\eta_4
    \right),
  \end{aligned}
\]
where $c_0 = a_0/2$, $c_1=a_1/3$ and $c_2=a_2/3$. With $\bm{v}(\eta_2) =
\bm{v}_1+\eta_2\bm{v}_2$ the first integral reads
\[
  \begin{aligned}
    J_1 &= \int\limits_0^1\bm{v}(\eta_2)\cdot\bm{n} \int\limits_{{(0,1)}^2}
    \frac{\eta_3 (c_0+c_1\eta_3+c_2\eta_3\eta_4)}
    {\seminorm{\eta_3\bm{u}_1+\eta_3\eta_4 \bm{u}_2 - \bm{v}(\eta_2)}^3}
    \,d\eta_3\,d\eta_4\,d\eta_2,
  \end{aligned}
\]
where the inner integral has the same form as the integrals of the previous
section. Therefore, we use~\eqref{eq:dlp_intedge} to compute the inner integral
and approximate the outer integral with numerical quadrature. We summarise the
parameters in Tables~\ref{table:dlp_parvert_p},~\ref{table:dlp_parvert_q}
and~\ref{table:dlp_parvert_d}. Note that for $J_2$ only $\bm{u}_i$ $\bm{v}_i$
need to be exchanged.

\subsubsection{Far-field}\label{subsubsec:dlp_far}
With the notation of Section~\ref{subsubsec:slp_far}, the integral on
${(0,1)}^4$ reads
\[
  J = \frac{g_\tau g_\sigma}{4\pi}\int\limits_{{(0,1)}^4} \eta_1\eta_3
  (a_0+a_1\eta_3+a_2\eta_3\eta_4)
  \frac{(-\bm{p} + \eta_1 \bm{v}_1- \eta_1 \eta_2 \bm{v}_2) \cdot \bm{n}}
  {\seminorm{\bm{p}+\eta_3 \bm{u}_1 + \eta_3 \eta_4\bm{u}_2 - \eta_1 \bm{v}_1
      - \eta_1 \eta_2 \bm{v}_2}^3}\,
  d\bm{\eta}.
\]
We integrate analytically with respect to $\eta_1$ and $\eta_2$ and 
numerically with respect to $\eta_3$ and $\eta_4$, i.e. we approximate $J$ by
\[
  \sum_{k,\ell=1}^n \omega_k\, \omega_\ell\, \eta^{(k)}
  (a_0+a_1\eta^{(k)}+a_2\eta^{(k)}\eta^{(\ell)})
  \int\limits_{{(0,1)}^2} 
  \frac{\eta_1 (-\bm{p} + \eta_1 \bm{v}_1- \eta_1 \eta_2 \bm{v}_2) \cdot \bm{n}}
  {\seminorm{\bm{u}(\eta^{(k)},\eta^{(\ell)}) - \eta_1 \bm{v}_1 -
      \eta_1 \eta_2 \bm{v}_2}^3}\,d\eta_1\,d\eta_2,
\]
with $\bm{u}(\eta^{(k)},\eta^{(\ell)}) = \bm{p}+\eta^{(k)} \bm{u}_1 +
\eta^{(k)} \eta^{(\ell)}\bm{u}_2$.
The parameters are listed in
Tables~\ref{table:dlp_parfar_p},~\ref{table:dlp_parfar_q}
and~\ref{table:dlp_parfar_d}.

\section{Numerical experiments}\label{sec:exp}
\begin{figure}[hbt]
  \centering
  \includegraphics[width=\linewidth]{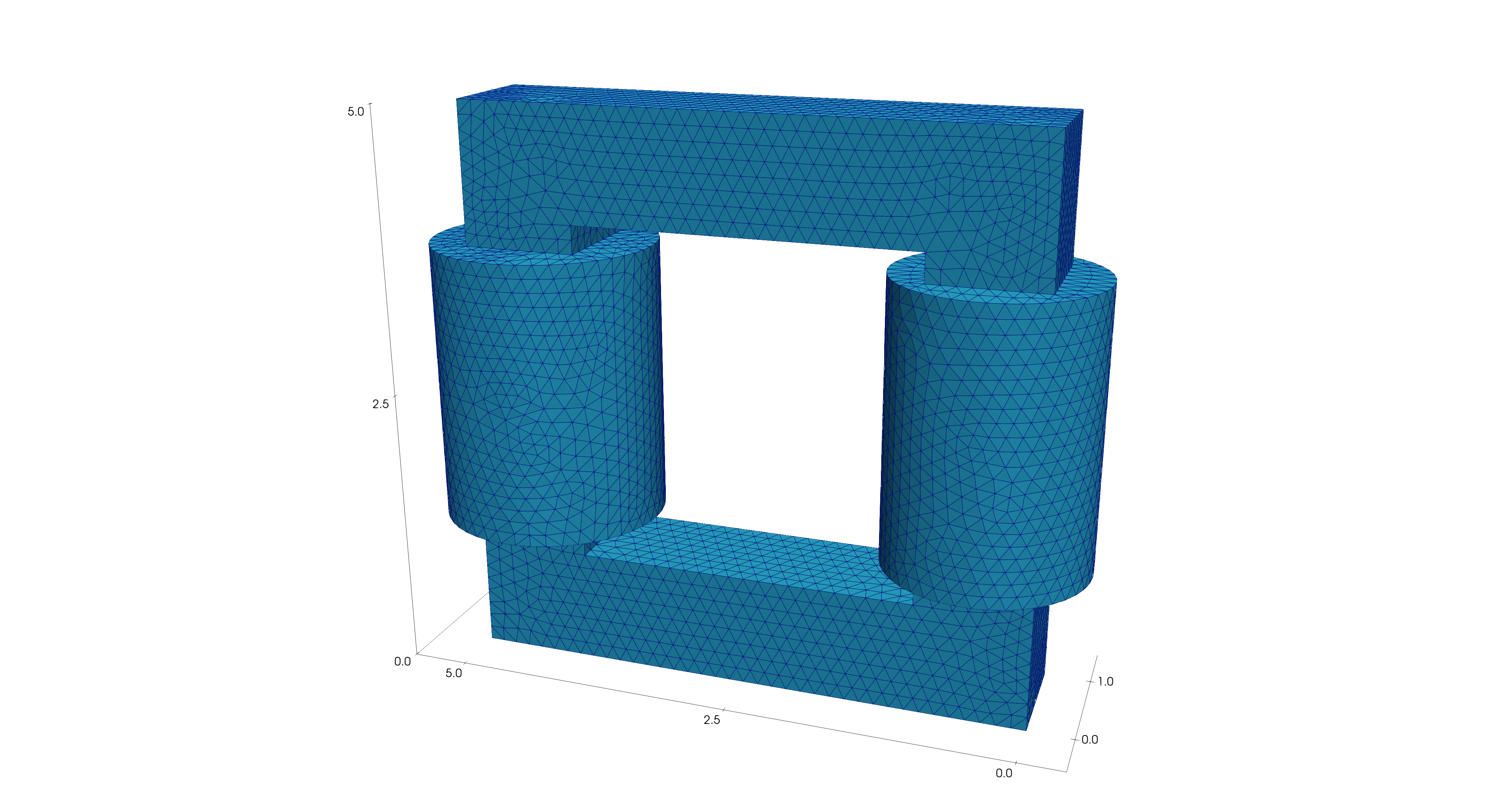}
  \caption{The surface $\Gamma_h^{(2)}$ consists of $10112$ triangles and
    resembles a transformer.
  }\label{fig:transformer}
\end{figure}

In this final section, we verify the correctness of the analytical formulae
in numerical examples. To this end, we consider two different geometries
for $\Gamma_h$, namely a triangulated unit sphere $\Gamma_h^{(1)}$ with $N=4608$
triangles and the surface $\Gamma_h^{(2)}$ of a transformer visualised in
Figure~\ref{fig:transformer}.

By $\bm{A}_{r,s}$ we denote the approximation of
the exact boundary element matrix $\bm{A}=\bm{V},\bm{K}$ computed by the
semi-analytical formulae. We use tensorised Gauss-Legendre quadrature rules with
$r\times r$ points for the far-field and $s$ points in the singular
vertex case. Hence, the computation of $\bm{V}_{r,s}$ requires
$\mathcal{O}(N s + N^2 r^2)$ evaluations of anti-derivatives. Note that the
conventional approach based on four-dimensional quadrature involves
$\mathcal{O}(N s^4 + N^2 r^4)$ kernel evaluations.

We measure the relative error
\[
  e = \norm{\bm{A}_{r,s} - \bm{A}}_{F} / \norm{\bm{A}}_{F}
\]
in the Frobenius norm defined by
\[
  \norm{\bm{A}}_F^2 = \sum_{m=1}^M \sum_{n=1}^N {\left(\bm{A}[m,n]\right)}^2,
  \quad \bm{A} \in \mathbb{R}^{M \times N}.
\]
Since the exact boundary element matrix $\bm{A}$ is not available, we compute a
reference approximation with four-dimensional quadrature of order $r=22$ and
$s=24$. Moreover, we set $s=r+2$ in all experiments.

\begin{figure}[hbt]
  \centering
  \includegraphics{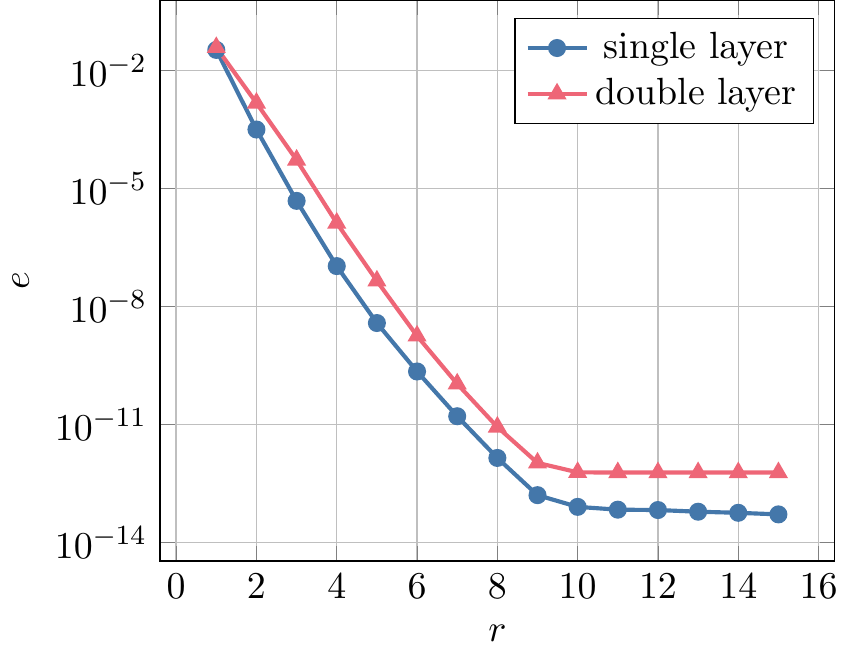}
  \caption{Relative error $e$ for increasing quadrature order $r$ for
    the sphere $\Gamma_h^{(1)}$.}\label{fig:frobsphere}
\end{figure}

\begin{figure}[hbt]
  \centering
  \includegraphics{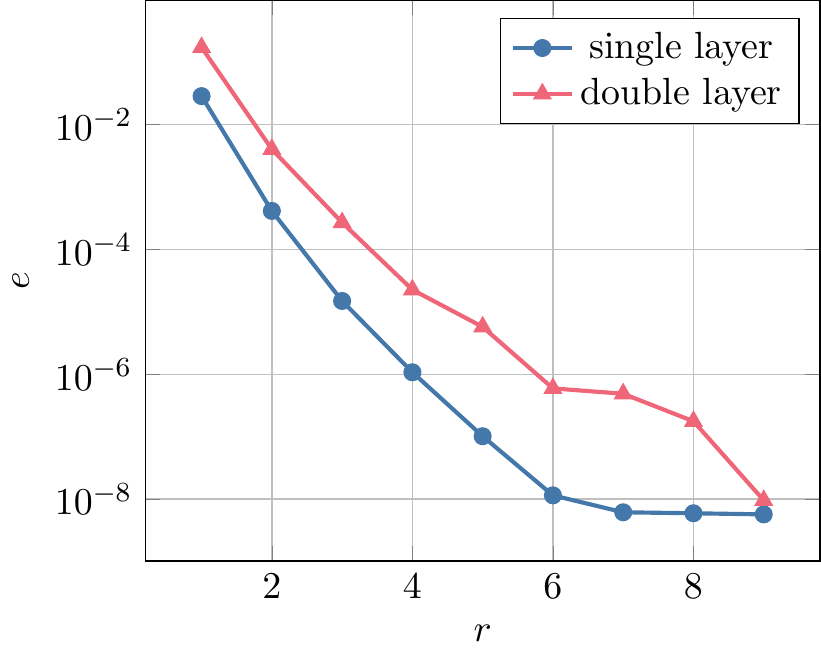}
  \caption{Relative error $e$ for increasing quadrature order $r$ for
    the transformer $\Gamma_h^{(2)}$.}\label{fig:frobtransformer}
\end{figure}

Figures~\ref{fig:frobsphere} and~\ref{fig:frobtransformer} show that the error
$e$ decreases exponentially in the quadrature order $r$ as expected. It reaches
$10^{-12}$ and $10^{-8}$ respectively and we see that $\bm{V}_{r,s}$ is slightly
more accurate than $\bm{K}_{r,s}$ for identical $r$ and $s$. Overall, we conclude
that the semi-analytical formulae produce numerically correct results.

\section{Conclusion}\label{sec:conc}
In comparison to black-box numerical quadrature, analytical integration exploits
the specific structure of the discretisation to reduce the computational costs
while preserving the level of accuracy. Since the regularisation method is not
limited to the Laplace equation, it is promising to extend the strategy to other
relevant kernel functions. Whereas the application to linear elasticity should
follow directly from the results presented here, the situation of time-harmonic
wave problems is less straight-forward due to oscillatory integrands. Whether
analytical integration is possible there needs to be investigated in future work.

\printbibliography%

@book{Prudnikov1981,
  AUTHOR =	 {Prudnikov, A. P. and Brychkov, Y. A. and Marichev, O. I.},
  TITLE =	 {Интегралы и ряды. Том 1. Integrals and series. Vol. 1.},
  PUBLISHER =	 {Fiziko-Matematicheskaya Literatura, Moscow},
  LANGUAGE =	 {langrussian},
  YEAR =	 2002,
  PAGES =	 632,
  ISBN =	 {5-9221-0323-7},
  MRCLASS =	 {00A22 (26-00 33-00 40-01)},
  MRNUMBER =	 2123874,
}

@book{Gradstein2015,
  AUTHOR =	 {Grad\v{s}te\u{\i}n, I. S. and Ry\v{z}ik, I. M.},
  TITLE =	 {Table of integrals, series, and products},
  EDITION =	 {Eighth},
  NOTE =	 {Translated from the Russian},
  PUBLISHER =	 {Elsevier/Academic Press, Amsterdam},
  YEAR =	 2015,
  PAGES =	 {xlvi + 1133},
  ISBN =	 {978-0-12-384933-5},
  MRCLASS =	 {00A22 (33-00)},
  MRNUMBER =	 3307944,
}

@book{McLean2000,
  AUTHOR =	 {McLean, W.},
  TITLE =	 {Strongly elliptic systems and boundary integral equations},
  PUBLISHER =	 {Cambridge University Press, Cambridge},
  YEAR =	 2000,
  PAGES =	 {xiv+357},
  ISBN =	 {0-521-66332-6; 0-521-66375-X},
  MRCLASS =	 {35J45 (47F05 47G10 47N20 65N38)},
  MRNUMBER =	 1742312,
  MRREVIEWER =	 {Dorina I. Mitrea},
}

@book{SauterSchwab2011,
  AUTHOR =	 {Sauter, S. A. and Schwab, C.},
  TITLE =	 {Boundary element methods},
  SERIES =	 {Springer Series in Computational Mathematics},
  VOLUME =	 39,
  NOTE =	 {Translated and expanded from the 2004 German original},
  PUBLISHER =	 {Springer-Verlag, Berlin},
  YEAR =	 2011,
  PAGES =	 {xviii+561},
  ISBN =	 {978-3-540-68092-5},
  MRCLASS =	 {65-02 (35J25 46N40 47G10 65N38)},
  MRNUMBER =	 2743235,
  MRREVIEWER =	 {Paul Andrew Martin},
  DOI =		 {10.1007/978-3-540-68093-2},
  URL =		 {https://doi.org/10.1007/978-3-540-68093-2},
}

@incollection{Sauter1998,
    AUTHOR = {Erichsen, S. and Sauter, S. A.},
     TITLE = {Efficient automatic quadrature in {$3$}-d {G}alerkin {BEM}},
      NOTE = {Seventh Conference on Numerical Methods and Computational
              Mechanics in Science and Engineering (NMCM 96) (Miskolc)},
   JOURNAL = {Comput. Methods Appl. Mech. Engrg.},
  FJOURNAL = {Computer Methods in Applied Mechanics and Engineering},
    VOLUME = {157},
      YEAR = {1998},
    NUMBER = {3-4},
     PAGES = {215--224},
      ISSN = {0045-7825},
   MRCLASS = {65N38},
  MRNUMBER = {1634288},
       DOI = {10.1016/S0045-7825(97)00236-3},
       URL = {https://doi.org/10.1016/S0045-7825(97)00236-3},
}

@article{Duffy1982,
    AUTHOR = {Duffy, M. G.},
     TITLE = {Quadrature over a pyramid or cube of integrands with a
              singularity at a vertex},
   JOURNAL = {SIAM J. Numer. Anal.},
  FJOURNAL = {SIAM Journal on Numerical Analysis},
    VOLUME = {19},
      YEAR = {1982},
    NUMBER = {6},
     PAGES = {1260--1262},
      ISSN = {0036-1429},
   MRCLASS = {65D30},
  MRNUMBER = {679664},
       DOI = {10.1137/0719090},
       URL = {https://doi.org/10.1137/0719090},
}

@book{FichtenholzII2003,
  TITLE =	 {Курс дифференциального и интегрального исчисления. Том
                  2. Differential and integral calculus. Vol. 2.},
  PUBLISHER =	 {Fiziko-Matematiches\-kaya Literatura, Moscow},
  author =	 {Fichtenholz, G. M.},
  LANGUAGE =	 {langrussian},
  EDITION =	 {Eighth},
  YEAR =	 2003,
  PAGES =	 863,
  ISBN =	 {5-9221-0155-9},
}

@book{Bebendorf2008,
  AUTHOR =	 {Bebendorf, M.},
  TITLE =	 {Hierarchical matrices},
  SERIES =	 {Lecture Notes in Computational Science and Engineering},
  VOLUME =	 63,
  PUBLISHER =	 {Springer, Berlin},
  YEAR =	 2008,
  PAGES =	 {xvi+290},
  ISBN =	 {978-3-540-77146-3},
  MRCLASS =	 {15-02 (35J40 65N22)},
  MRNUMBER =	 2451321,
  MRREVIEWER =	 {Elena Pelican},
}

@book {RjasanowSteinbach2007,
    AUTHOR = {Rjasanow, S. and Steinbach, O.},
     TITLE = {The fast solution of boundary integral equations},
    SERIES = {Mathematical and Analytical Techniques with Applications to
              Engineering},
 PUBLISHER = {Springer, New York},
      YEAR = {2007},
     PAGES = {xii+279},
      ISBN = {978-0-387-34041-8; 0-387-34041-6},
   MRCLASS = {65N38},
  MRNUMBER = {2310663},
MRREVIEWER = {Mar\'{\i}a-Luisa Rap\'{u}n},
}

@article{SchwabWendland1992,
    AUTHOR = {Schwab, C. and Wendland, W. L.},
     TITLE = {On numerical cubatures of singular surface integrals in
              boundary element methods},
   JOURNAL = {Numer. Math.},
  FJOURNAL = {Numerische Mathematik},
    VOLUME = {62},
      YEAR = {1992},
    NUMBER = {3},
     PAGES = {343--369},
      ISSN = {0029-599X},
   MRCLASS = {65D30 (45D05 47G30 65N38)},
  MRNUMBER = {1169009},
MRREVIEWER = {Stephen W. Brady},
}

@article {MousaviSukumar2010,
    AUTHOR = {Mousavi, S. E. and Sukumar, N.},
     TITLE = {Generalized {D}uffy transformation for integrating vertex
              singularities},
   JOURNAL = {Comput. Mech.},
  FJOURNAL = {Computational Mechanics},
    VOLUME = {45},
      YEAR = {2010},
    NUMBER = {2-3},
     PAGES = {127--140},
      ISSN = {0178-7675},
   MRCLASS = {65D30 (65N30 65N38)},
  MRNUMBER = {2576090},
       DOI = {10.1007/s00466-009-0424-1},
       URL = {https://doi.org/10.1007/s00466-009-0424-1},
}

@article {Taylor2003,
    AUTHOR = {Taylor, D. J.},
     TITLE = {Accurate and efficient numerical integration of weakly
              singular integrals in {G}alerkin {EFIE} solutions},
   JOURNAL = {IEEE Trans. Antennas and Propagation},
  FJOURNAL = {Institute of Electrical and Electronics Engineers.
              Transactions on Antennas and Propagation},
    VOLUME = {51},
      YEAR = {2003},
    NUMBER = {7},
     PAGES = {1630--1637},
      ISSN = {0018-926X},
   MRCLASS = {78M25 (65D30)},
  MRNUMBER = {2004244},
       DOI = {10.1109/TAP.2003.813623},
       URL = {https://doi.org/10.1109/TAP.2003.813623},
}

@book{DautrayLions1990,
    AUTHOR = {Dautray, R. and Lions, J.-L.},
     TITLE = {Mathematical analysis and numerical methods for science and
              technology. {V}ol. 3},
 PUBLISHER = {Springer-Verlag, Berlin},
      YEAR = {1990},
     PAGES = {x+515},
      ISBN = {3-540-50208-4; 3-540-66099-2},
   MRCLASS = {00A05 (35-01 47-01 78-01 81-01)},
  MRNUMBER = {1064315},
}
\appendix
\section{Appendix}
In the following, the parameter tables for the formulae of the double layer
potential are listed.
\vfill
\begin{table}[htb]
  \centering
  \begin{tabular}{c|c||c|c|c}
    $i$ & $j$ & $q$ & $q\cdot q_0$ & $2q\cdot q_1$\\
    \hline\hline
    $1$ & $0$ & $\seminorm{\bm{u}+\bm{v}}^2$ & $\seminorm{\bm{u}}^2$
    & $-(\bm{u}+\bm{v})\cdot\bm{u}$\\
    \hline
    & $1$ & $\seminorm{\bm{u}+\bm{v}}^2$ & $\seminorm{\bm{u}-\bm{w}}^2$
    & $-(\bm{u}+\bm{v})\cdot(\bm{u}-\bm{w})$\\
    \hline
    $2$ & $0$ & $\seminorm{\bm{u}}^2$ & $\seminorm{\bm{w}}^2$
    & $-\bm{u}\cdot \bm{w}$\\
    \hline
    & $1$ &$\seminorm{\bm{v}}^2$ & $\seminorm{\bm{w}}^2$
    & $\bm{v}\cdot \bm{w}$\\
    \hline
    $3$ & $0$ & $\seminorm{\bm{w}}^2$ &
    $\seminorm{\bm{u}}^2$ & $-\bm{u}\cdot\bm{w}$ \\
    \hline
    & $1$ & $\seminorm{\bm{v}}^2$ & $\seminorm{\bm{u}}^2$ & $\bm{u}\cdot\bm{v}$ \\
    \hline
    $4$ & $0$ & $\seminorm{\bm{v}+\bm{w}}^2$ & $\seminorm{\bm{w}}^2$
    &$-(\bm{v}+\bm{w})\cdot\bm{w}$\\
    \hline
    & $1$ & $\seminorm{\bm{u}+\bm{v}+\bm{w}}^2$ & $\seminorm{\bm{w}}^2$
    &$-(\bm{u}+\bm{v}+\bm{w})\cdot\bm{w}$\\
    \hline
    $5$ & $0$ & $\seminorm{\bm{u}}^2$ & $\seminorm{\bm{w}}^2$
    &$-\bm{u}\cdot\bm{w}$\\
    \hline
    & $1$ & $\seminorm{\bm{u}+\bm{v}+\bm{w}}^2$ & $\seminorm{\bm{w}}^2$
    &$-(\bm{u}+\bm{v}+\bm{w})\cdot\bm{w}$\\
    \hline 
  \end{tabular}
  \caption{Parameters $q$ and $q_k$ in~\eqref{eq:dlp_intedge} for the integrals
    of $h^{(i)}_j$ for the edge case.}\label{table:dlp_paredge_q}
\end{table}
\vfill
\begin{table}[htb]
  \centering
  \begin{tabular}{c||c|c|c}
    $i$&$d$&$d\cdot d_0$&$2d\cdot d_1$\\
    \hline\hline
    $1$
    &$\begin{aligned}
      &\seminorm{\bm{u}+\bm{v}}^2\seminorm{\bm{w}}^2\\
      &-{((\bm{u}+\bm{v})\cdot\bm{w})}^2
    \end{aligned}$
    &$\begin{aligned}
      &\seminorm{\bm{u}}^2\seminorm{\bm{w}}^2-{(\bm{u}\cdot\bm{w})}^2
    \end{aligned}$
    &$\begin{aligned}
      &{(\bm{u}\cdot\bm{w})}\, (\bm{u}+\bm{v})\cdot\bm{w}\\
      &-\seminorm{\bm{w}}^2(\bm{u}+\bm{v})\cdot\bm{u}
    \end{aligned}$\\
    \hline
    $2$
    &$\begin{aligned}
      &\seminorm{\bm{u}}^2\seminorm{\bm{v}}^2-{(\bm{u}\cdot\bm{v})}^2
    \end{aligned}$
    &$\begin{aligned}
      &\seminorm{\bm{u}+\bm{v}}^2\seminorm{\bm{w}}^2\\
      &-{((\bm{u}+\bm{v})\cdot\bm{w})}^2
    \end{aligned}$
    &$\begin{aligned}
      &\bm{v}\cdot\bm{w}\, (\bm{u}+\bm{v})\cdot\bm{u}\\
      &-\bm{u}\cdot\bm{w}\ (\bm{u}+\bm{v})\cdot\bm{v}
    \end{aligned}$\\
    \hline
    $3$
    &$\begin{aligned}
      &\seminorm{\bm{v}}^2\seminorm{\bm{w}}^2-{(\bm{v}\cdot\bm{w})}^2
    \end{aligned}$
    &$\begin{aligned}
      &\seminorm{\bm{v}+\bm{w}}^2\seminorm{\bm{u}}^2\\
      &-{((\bm{v}+\bm{w})\cdot\bm{u})}^2
    \end{aligned}$
    &$\begin{aligned}
      &\bm{u}\cdot\bm{v}\ (\bm{v}+\bm{w})\cdot\bm{w}\\
      &-\bm{u}\cdot\bm{w}\ (\bm{v}+\bm{w})\cdot\bm{v}
    \end{aligned}$\\
    \hline
    $4$
    &$\begin{aligned}
      &\seminorm{\bm{u}}^2\seminorm{\bm{v}+\bm{w}}^2\\
      &-{((\bm{v}+\bm{w})\cdot\bm{u})}^2
    \end{aligned}$
    &$\begin{aligned}
      &\seminorm{\bm{u}}^2\seminorm{\bm{w}}^2-{(\bm{u}\cdot\bm{w})}^2
    \end{aligned}$
    &$\begin{aligned}
      &\bm{u}\cdot\bm{w}\ (\bm{v}+\bm{w})\cdot\bm{u}\\
      &-\seminorm{\bm{u}}^2 (\bm{v}+\bm{w})\cdot\bm{w}
    \end{aligned}$\\
    \hline
    $5$
    &$\begin{aligned}
      &\seminorm{\bm{u}}^2\seminorm{\bm{v}+\bm{w}}^2\\
      &-{((\bm{v}+\bm{w})\cdot\bm{u})}^2
    \end{aligned}$
    &$\begin{aligned}
      &\seminorm{\bm{v}}^2\seminorm{\bm{w}}^2-{(\bm{v}\cdot\bm{w})}^2
    \end{aligned}$
    &$\begin{aligned}
      &\bm{u}\cdot\bm{v}\ (\bm{v}+\bm{w})\cdot\bm{w}\\
      &-\bm{u}\cdot\bm{w}\ (\bm{v}+\bm{w})\cdot\bm{v}
    \end{aligned}$\\
    \hline 
  \end{tabular}
  \caption{Parameters $d$ and $d_k$ in~\eqref{eq:dlp_intedge} for the integrals
    of $h^{(i)}_j$ for the edge case.}\label{table:dlp_paredge_d}
\end{table}

\begin{sidewaystable}[htb]
  \centering
  \begin{tabular}{c|c||c|c|c|c}
    $i$ & $j$ &$p_0$ & $p_1$ & $p_2$ & $p_3$\\
    \hline\hline
    $1$ & $0$ & $-(c_0+c_2)\seminorm{\bm{u}}^2$
    &$\begin{aligned}
      &(2c_0-c_1+3c_2)\seminorm{\bm{u}}^2\\
      &+2(c_0+c_2)\bm{u}\cdot\bm{v}
    \end{aligned}$
    &$\begin{aligned}
      &(2c_1-c_0-3c_2)\seminorm{\bm{u}}^2-(c_0+c_2)\seminorm{\bm{v}}^2\\
      &+2(c_1-c_0-2c_2)\bm{u}\cdot\bm{v}
    \end{aligned}$
    &$(c_2-c_1)\seminorm{\bm{u}+\bm{v}}^2$\\
    \hline
    & $1$ & $(c_0+c_2) (\bm{w}-\bm{u})\cdot \bm{u}$
    &$\begin{aligned}
      &2(c_0+c_2)\bm{u}\cdot\bm{v}-(c_0+c_2)\bm{v}\cdot\bm{w}\\
      &+(2c_0-c_1+3c_2)\seminorm{\bm{u}}^2\\
      &+(c_1-c_0-2c_2)\bm{u}\cdot\bm{w}
    \end{aligned}$
    &$\begin{aligned}
      &(2c_1-c_0-3c_2)\seminorm{\bm{u}}^2-(c_0+c_2)\seminorm{\bm{v}}^2\\
      &+2(c_1-c_0-2c_2)\bm{u}\cdot\bm{v}\\
      &+(c_2-c_1)(\bm{u}+\bm{v})\cdot\bm{w}
    \end{aligned}$
    & $(c_2-c_1)\seminorm{\bm{u}+\bm{v}}^2$\\
    \hline
    $2$ & $0$
    &$\begin{aligned}
      &c_0(\bm{u}+\bm{v})\cdot\bm{w}\\
      &+(c_2-c_1)\seminorm{\bm{w}}^2
    \end{aligned}$
    &$\begin{aligned}
      &(2c_1-c_2)\bm{u}\cdot\bm{w}+c_2\bm{v}\cdot\bm{w}\\
      &-c_0(\bm{u}+\bm{v})\cdot\bm{u}      
    \end{aligned}$
    & $-(c_1\seminorm{\bm{u}}^2+c_2\bm{u}\cdot\bm{v})$
    &$0$\\
    \hline
    & $1$
    &$\begin{aligned}
      &c_0(\bm{u}+\bm{v})\cdot\bm{w}\\
      &+(c_2-c_1)\seminorm{\bm{w}}^2
    \end{aligned}$
    &$\begin{aligned}
      &(2c_2-c_1)\bm{v}\cdot\bm{w}+c_1\bm{u}\cdot\bm{w}\\
      &+c_0(\bm{u}+\bm{v})\cdot\bm{v}
    \end{aligned}$
    & $c_2\seminorm{\bm{v}}^2+c_1\bm{u}\cdot\bm{v}$
    &$0$\\
    \hline
    $3$ & $0$
    &$(c_0+c_2) \seminorm{\bm{u}}^2$
    &$(c_0+c_2) (\bm{v}-\bm{w})\cdot\bm{u}$
    &$-(c_0+c_2) \bm{v}\cdot\bm{w}$
    &$0$\\
    \hline
    & $1$
    &$(c_0+c_2) \seminorm{\bm{u}}^2$
    &$2(c_0+c_2)\bm{u}\cdot\bm{v}$
    &$(c_0+c_2) \seminorm{\bm{v}}^2$
    &$0$\\
    \hline
    $4$ & $0$ &$-c_0\bm{u}\cdot\bm{w}-c_2\seminorm{\bm{w}}^2$
    & $\begin{aligned}
      &c_0(\bm{v}+2\bm{w})\cdot\bm{u}\\
      &+c_2 (2\bm{v}+3\bm{w})\cdot\bm{w}
    \end{aligned}$
    & $\begin{aligned}
      &-c_0(\bm{v}+\bm{w})\cdot\bm{u}\\
      &- c_2(\seminorm{\bm{v}+\bm{w}}^2+2(\bm{v}+\bm{w})\cdot\bm{w}) 
    \end{aligned}$
    & $c_2\seminorm{\bm{v}+\bm{w}}^2$\\
    \hline
    & $1$ &$-c_0\bm{u}\cdot\bm{w}-c_2\seminorm{\bm{w}}^2$
    & $\begin{aligned}
      &c_0 (\bm{u}+\bm{v})\cdot\bm{u}+(2c_0+c_2)\bm{u}\cdot\bm{w}\\
      &+c_2 (2\bm{v}+3\bm{w})\cdot\bm{w}
      \end{aligned}$
    & $\begin{aligned}
      &-c_0(\bm{u}+\bm{v}+\bm{w})\cdot\bm{u}-c_2\bm{u}\cdot\bm{v}\\
      &-c_2\seminorm{\bm{v}+\bm{w}}^2 -2c_2(\bm{u}+\bm{v}+\bm{w})\cdot\bm{w}
    \end{aligned}$
    & $\begin{aligned}
      &c_2(\bm{v}+\bm{w})\cdot\bm{u}\\
      &+c_2\seminorm{\bm{v}+\bm{w}}^2
    \end{aligned}$\\
    \hline
    $5$ & $0$ & $-c_0\bm{v}\cdot\bm{w}$
    &$c_0(\bm{v}-\bm{w})\cdot\bm{u}-c_2\bm{v}\cdot\bm{w}$
    &$c_0\seminorm{\bm{u}}^2+c_2(\bm{v}-\bm{w})\cdot\bm{u}$
    &$c_2\seminorm{\bm{u}}^2$\\
    \hline
    & $1$ & $-c_0\bm{v}\cdot\bm{w}$
    &$\begin{aligned}
      &c_0(\bm{u}+\bm{v})\cdot\bm{v}-c_2\bm{v}\cdot\bm{w}\\
      &+c_0(\bm{v}-\bm{u})\cdot\bm{w}
    \end{aligned}$
    &$\begin{aligned}
      &c_0(\bm{u}+\bm{v}+\bm{w})\cdot\bm{u}\\
      &+c_2\left((\bm{u}+\bm{v})\cdot\bm{v}+(\bm{v}-\bm{u})\cdot\bm{w}\right)
    \end{aligned}$
    &$c_2(\bm{u}+\bm{v}+\bm{w})\cdot\bm{u}$\\
    \hline 
  \end{tabular}
  \caption{Parameters $p_k$ in~\eqref{eq:dlp_intedge} for the integrals
    of $h^{(i)}_j$ for the edge case.}\label{table:dlp_paredge_p}
\end{sidewaystable}


\begin{sidewaystable}[htb]
  \centering
  \begin{tabular}{c|c||c|c|c|c}
    $i$ & $j$ &$p_0$ & $p_1$ & $p_2$ & $p_3$\\
    \hline\hline
    $1$ & $0$
    &$-c_0 \bm{u}_2\cdot\bm{v} - c_2 \seminorm{\bm{v}}^2$
    &$c_0 \bm{u}_1\cdot\bm{u}_2+2 c_2 \bm{u}_1\cdot\bm{v} - c_1 \bm{u}_2\cdot\bm{v}$
    &$-c_2 \seminorm{\bm{u}_1}^2+c_1 \bm{u}_1\cdot\bm{u}_2$
    &$0$\\
    \hline
    &$1$
    &$-c_0 \bm{u}_2\cdot\bm{v}-c_2 \seminorm{\bm{v}}^2$
    &$2 c_2 \bm{u}_1\cdot\bm{v}+c_0 \bm{u}_2\cdot(\bm{u}_1+\bm{u}_2)+(c_2-c_1) \bm{u}_2\cdot\bm{v}$
    &$-c_2 \bm{u}_1\cdot (\bm{u}_1+\bm{u}_2)+c_1 \bm{u}_2\cdot(\bm{u}_1+\bm{u}_2)$
    $0$\\
    \hline
  \end{tabular}
  \caption{Parameters $p_k$ in~\eqref{eq:dlp_intedge} for the integrals
    of $h^{(i)}_j$ for the vertex case.}\label{table:dlp_parvert_p}
  \vspace{3cm}
  \begin{tabular}{c|c||c|c|c|c}
    $i$ & $j$ &$p_0$ & $p_1$ & $p_2$ & $p_3$\\
    \hline\hline
    $1$ & $0$
    &$-\bm{n}\cdot \bm{v}_2 \seminorm{\bm{u}}^2 +
     (\bm{n}\cdot \bm{p})(\bm{u}\cdot \bm{v}_2)$
    &$\begin{aligned}
      &2(\bm{n}\cdot\bm{v}_2)(\bm{u}\cdot\bm{v}_1)-
      (\bm{n}\cdot\bm{v}_1)(\bm{u}\cdot\bm{v}_2)\\
      &-(\bm{n}\cdot\bm{p})(\bm{v}_1\cdot\bm{v}_1)
    \end{aligned}$
    &$\begin{aligned}
      -(\bm{n}\cdot\bm{v}_2)\seminorm{\bm{v}_1}^2+
      (\bm{n}\cdot\bm{v}_1)(\bm{v}_1\cdot\bm{v}_2)
    \end{aligned}$
    &$0$\\
    \hline
    & $1$
    &$-\bm{n}\cdot \bm{v}_2 \seminorm{\bm{u}}^2+
     (\bm{n}\cdot \bm{p})(\bm{u}\cdot \bm{v}_2)$
    &$\begin{aligned}
      &(\bm{n}\cdot\bm{v}_2)(2\bm{u}\cdot\bm{v}_1+\bm{u}\cdot\bm{v}_2)-
      (\bm{n}\cdot\bm{v}_1)(\bm{u}\cdot\bm{v}_2)\\
      &-(\bm{n}\cdot\bm{p})(\bm{v}_1\cdot\bm{v}_2+\seminorm{\bm{v}_2}^2)
    \end{aligned}$
    &$\begin{aligned}
      &-(\bm{n}\cdot\bm{v}_2)(\seminorm{\bm{v}_1}^2+\bm{v}_1\cdot\bm{v}_2)\\
      &+(\bm{n}\cdot\bm{v}_1)(\seminorm{\bm{v}_2}^2+\bm{v}_1\cdot\bm{v}_2)
    \end{aligned}$
    &$0$\\
    \hline
  \end{tabular}
  \caption{Parameters $p_k$ in~\eqref{eq:dlp_intedge} for the integrals
    of $h^{(i)}_j$ for the far-field case.}\label{table:dlp_parfar_p}
\end{sidewaystable}


\begin{table}[htb]
  \centering
  \begin{tabular}{c|c||c|c|c}
    $i$ & $j$ & $q$ & $q\cdot q_0$ & $2q\cdot q_1$\\
    \hline\hline
    $1$ & $0$ & $\seminorm{\bm{u}_1}^2$ & $\seminorm{\bm{v}}^2$ & $\bm{u}_1\cdot\bm{v}$
    \\
    \hline
    & $1$ & $\seminorm{\bm{u}_1+\bm{u}_2}^2$ & $\seminorm{\bm{v}}^2$ & $-(\bm{u}_1+\bm{u}_2)\cdot\bm{v}$
    \\
  \end{tabular}
  \caption{Parameters $q$ and $q_k$ in~\eqref{eq:dlp_intedge} for the integrals
    of $h^{(i)}_j$ for the vertex case.}\label{table:dlp_parvert_q}
\end{table}

\begin{table}[htb]
  \centering
  \begin{tabular}{c||c|c|c}
    $i$&$d$&$d\cdot d_0$&$2d\cdot d_1$\\
    \hline\hline
    $1$
    &$\seminorm{\bm{u}_1}^2\seminorm{\bm{u}_2}^2-{(\bm{u}_1\cdot\bm{u}_2)}^2$
    &$\seminorm{\bm{u}_2}^2\seminorm{\bm{v}}^2-{(\bm{u}_2\cdot\bm{v})}^2$
    &$(\bm{u}_1\cdot\bm{u}_2)(\bm{u}_2\cdot\bm{v})-\seminorm{\bm{u}_2}^2\bm{u}_1\cdot\bm{v}$\\
    \hline
  \end{tabular}
  \caption{Parameters $d$ and $d_k$ in~\eqref{eq:dlp_intedge} for the integrals
    of $h^{(i)}_j$ for the vertex case.}\label{table:dlp_parvert_d}
\end{table}

\begin{table}[htb]
  \centering
  \begin{tabular}{c|c||c|c|c}
    $i$ & $j$ & $q$ & $q\cdot q_0$ & $2q\cdot q_1$\\
    \hline\hline
    $1$ & $0$ & $\seminorm{\bm{u}+\bm{v}}^2$ & $\seminorm{\bm{u}}^2$
    & $-(\bm{u}+\bm{v})\cdot\bm{u}$\\
    \hline
    & $1$ & $\seminorm{\bm{u}+\bm{v}}^2$ & $\seminorm{\bm{u}-\bm{w}}^2$
    & $-(\bm{u}+\bm{v})\cdot(\bm{u}-\bm{w})$\\
    \hline
  \end{tabular}
  \caption{Parameters $q$ and $q_k$ in~\eqref{eq:dlp_intedge} for the integrals
    of $h^{(i)}_j$ for the far-field case.}\label{table:dlp_parfar_q}
\end{table}

\begin{table}[htb]
  \centering
  \begin{tabular}{c||c|c|c}
    $i$&$d$&$d\cdot d_0$&$2d\cdot d_1$\\
    \hline\hline
    $1$
    &$\begin{aligned}
      &\seminorm{\bm{u}+\bm{v}}^2\seminorm{\bm{w}}^2\\
      &-{((\bm{u}+\bm{v})\cdot\bm{w})}^2
    \end{aligned}$
    &$\begin{aligned}
      &\seminorm{\bm{u}}^2\seminorm{\bm{w}}^2-{(\bm{u}\cdot\bm{w})}^2
    \end{aligned}$
    &$\begin{aligned}
      &{(\bm{u}\cdot\bm{w})}\, (\bm{u}+\bm{v})\cdot\bm{w}\\
      &-\seminorm{\bm{w}}^2(\bm{u}+\bm{v})\cdot\bm{u}
    \end{aligned}$\\
    \hline
  \end{tabular}
  \caption{Parameters $d$ and $d_k$ in~\eqref{eq:dlp_intedge} for the integrals
    of $h^{(i)}_j$ for the far-field case.}\label{table:dlp_parfar_d}
\end{table}

\end{document}